\documentclass{amsart}
\usepackage{amsmath} 
\usepackage{amssymb}
\usepackage{mathrsfs}

\newtheorem{theorem}{Theorem}[section] 
\newtheorem{claim}[theorem]{Claim}

\newtheorem{conclusion}[theorem]{Conclusion}
\newtheorem{observation}[theorem]{Observation}

\theoremstyle{definition}
\newtheorem{definition}[theorem]{Definition}
\newtheorem{fact}[theorem]{Fact}
\newtheorem{convention}[theorem]{Convention}

\newtheorem{discussion}[theorem]{Discussion}

\newtheorem{hypothesis}[theorem]{Hypothesis}

\theoremstyle{remark}
\newtheorem{remark}[theorem]{Remark}
\newtheorem{question}[theorem]{Question}
\newtheorem{notation}[theorem]{Notation}

\newcommand{\Ker}{{\rm Ker}}
\newcommand{\Th}{{\rm Th}}
\newcommand{\tp}{{\rm tp}}
\newcommand{\st}{{\rm st}}

\newcommand{\ba}{{\rm ba}}

\newcommand{\fr}{{\rm fr}}
\newcommand{\Av}{{\rm Av}}
\newcommand{\acl}{{\rm acl}}
\newcommand{\voc}{{\rm voc}}
\newcommand{\stp}{{\rm stp}}

\newcommand{\spec}{{\rm spec}}

\newcommand{\uf}{{\rm uf}}

\newcommand{\iif}{{\rm if}}

\newcommand{\eq}{{\rm eq}}

\newcommand{\trufil}{{\rm tr-ufil}}

\newcommand{\EQ}{{\rm EQ}}

\newcommand{\JEP}{{\rm JEP}}

\newcommand{\Mod}{{\rm Mod}}
\newcommand{\Rang}{{\rm Rang}}

\newcommand{\rest}{{\restriction}}
\newcommand{\dom}{{\rm dom}}

\newcommand{\wilog}{{\rm without loss of generality}}
\newcommand{\Wilog}{{\rm Without loss of generality}}

\newcommand{\then}{{\underline{then}}}
\newcommand{\when}{{\underline{when}}}
\newcommand{\Then}{{\underline{Then}}}
\newcommand{\but}{{\underline{but}}}

\newcommand{\Iff}{{\underline{iff}}}
\newcommand{\mn}{{\medskip\noindent}}
\newcommand{\sn}{{\smallskip\noindent}}

\newcommand{\bfj}{{\bold j}}

\newcommand{\bfN}{{\bold N}}
\newcommand{\bfi}{{\bold i}}
\newcommand{\bft}{{\bold t}}

\newcommand{\bfS}{{\bold S}}

\newcommand{\bbB}{{\mathbb B}}
\newcommand{\bbT}{{\mathbb T}}

\newcommand{\gC}{{\mathfrak C}}

\newcommand{\cI}{{\mathscr I}}

\newcommand{\bbL}{{\mathbb L}}

\newcommand{\cM}{{\mathscr M}}

\newcommand{\cP}{{\mathscr P}}

\newcommand{\varp}{{\varepsilon}}

\newcommand{\cf}{{\rm cf}}

\newcount\skewfactor
\def\mathunderaccent#1#2 {\let\theaccent#1\skewfactor#2
\mathpalette\putaccentunder}
\def\putaccentunder#1#2{\oalign{$#1#2$\crcr\hidewidth
\vbox to.2ex{\hbox{$#1\skew\skewfactor\theaccent{}$}\vss}\hidewidth}}

\newenvironment{PROOF}[2][\proofname.]
   {\begin{proof}[#1]\renewcommand{\qedsymbol}{\textsquare$_{\rm #2}$}}
   {\end{proof}}

\begin{document}

\title {Hanf number for the strictly stable cases}
\author {Saharon Shelah}
\address{Einstein Institute of Mathematics\\
Edmond J. Safra Campus, Givat Ram\\
The Hebrew University of Jerusalem\\
Jerusalem, 91904, Israel\\
 and \\
 Department of Mathematics\\
 Hill Center - Busch Campus \\ 
 Rutgers, The State University of New Jersey \\
 110 Frelinghuysen Road \\
 Piscataway, NJ 08854-8019 USA}
\email{shelah@math.huji.ac.il}
\urladdr{http://shelah.logic.at}
\thanks{Partially supported by European Research Council Grant \#338821.
The author thanks Alice Leonhardt for the beautiful typing.
First typed December 28, 2012. Paper 1048}





\date{December 13, 2018}

\begin{abstract}
Suppose $\bold t = (T,T_1,p)$ is a triple of two first order 
theories $T \subseteq
T_1$ in vocabularies $\tau \subseteq \tau_1$ (respectively) 
of cardinality $\lambda$ and a
$\tau_1$-type $p$ over the empty set; the main case here is with $T$
stable.   We show that the Hanf number for the property: ``there is a 
model $M_1$ of $T_1$ which omits $p$, but $M_1
\rest \tau$ is saturated" is larger than the Hanf number of
$\bbL_{\lambda^+,\kappa}$ but smaller than the Hanf number of
$\bbL_{(2^\lambda)^+,\kappa}$ when $T$ is stable with $\kappa =
\kappa(T)$.  In fact, we characterize the Hanf number of $\bold t$
when we fix $(T,\lambda)$ where $T$ is a 
first order complete, $\lambda \ge |T|$ and demand $|T_1| \le \lambda$.
\end{abstract}

\maketitle
\numberwithin{equation}{section}
\setcounter{section}{-1}
\newpage

\section {Introduction}
\bigskip

\subsection {Background on Results}\
\bigskip

This continues papers of Baldwin-Shelah, starting from a problem of
Newelski \cite{New12} concerning the Hanf number
described in the abstract for classes $\bold t \in
\bold N_{\lambda,T}$ (defined formally in \ref{b2}), that is:
\mn
\begin{itemize}
\item  for $T$ is a complete first order theory, $\lambda$ an infinite
  cardinal $\ge |T|$ let $\bfN_{\lambda,T}$ be the class of triples
  $\bft = (T,T_1,p)$ such that $T_1 \supseteq T$ is first order of
  cardinality $\le \lambda$ and $p=p(x)$ a type in the vocabulary of $T_1$
\sn
\item  for $\bft \in \bfN_{\lambda,T},M$ is a model of $\bft$ \Iff
  \, it is a model of $T_1$ (so have the same vocabulary) omitting the
type $p$ such that its restriction to the vocabulary of $T$ is a
  saturated model
\sn
\item  the Hanf number $H(\bft)$ of $\bft \in \bfN_{\lambda,T}$ is the
  first cardinal $\mu$ such that $\bft$ has no model of cardinaity
  $\ge \mu$ and is infinity when there is no such bound
\sn
\item  the Hanf number $H(\bfN_{\lambda,T})$ of $\bfN_{\lambda,T}$ is
  $\sup\{H(\bft):\bft \in \bfN_{\lambda,T}$ and $H(\bft) < \infty\}$
\sn
\item  the Hanf number $H_{\bfN}(\lambda)$ is
  $\sup\{H(\bfN_{T,\lambda}):(T,\lambda)$ as above$\}$
\sn
\item  note that, considering $\bfN_{\lambda,T}$ if $T$ is unstable it
  is natural to assume that $\{\mu:\mu = \mu^{< \mu}\}$ is an
  unbounded class as otherwise for any $T_1,\lambda$ we have
  $H((T_1,T,\lambda)) \le \sup \{\mu^+:\mu = \mu^{< \mu}\}$;
Newelski in \cite{New12} essentially asks what is $H_{\bfN}(\lambda)$,
Baldwin-Shelah \cite{BlSh:958}, \cite{BlSh:992} have dealt with those numbers.
\end{itemize}
\mn
They showed in \cite{BlSh:958} that the Hanf number
$H_{\bfN}(\lambda)$ is essentially equal to the L\"owenheim number of
second order logic using unstable $T$'s and in \cite{BlSh:992} showed 
that for superstable $T,H(\bfN_{\lambda,T})$ is bigger than the 
Hanf number of $\bbL_{(2^\lambda)^+,\aleph_0}$
but it is smaller than $\bbL_{\beth_2(\lambda)^+,\aleph_0}$.

Our original aim was to deal with the case where $T$ is a stable theory and
concentrate on the strictly stable case (i.e. stable not superstable).  

However, we ask a stronger question.
\begin{question}
\label{x2}
Fix a complete first order theory $T$ and a cardinal $\lambda \ge
|T|$, what is $H(\bfN_{\lambda,T})$? recalling it is 
$\sup\{H(\bold t):H(\bold t) < \infty$ and 
$\bold t$ as above with $T_{\bold t} = T$ and $|T_{\bold t,1}| 
\le \lambda$, i.e. belongs to $\bold N_{\lambda,T}$ from
\ref{b2}(1)$\}$, recalling $H(\bold t)$ is the
supremum of the cardinalities of models in $\Mod_{\bold t}$.

Clearly this is a considerably more ambitious question.  
Now \cite{BlSh:958} actually determines $H(\bfN_{\lambda,T})$ when $T$
is unstable, so we shall concentrate here on the case $T$ is stable.
We give a quite complete answer.  For $T$
strictly stable, our original case, it appears that only the 
cardinals $|T|,\kappa(T)$ and a derived 
Boolean Algebra $\bbB(T)$ of cardinality $|D(T)|$, and a little more
where $D(T) = \cup\{D_n(T):n < \omega\},D_n(T)$ is the set of complete
$n$-types realized in models of $T$.
In fact, for any $T$, the little more is the truth value of 
$(2^{\aleph_0} > |D(T)| > |T| \wedge ``T$ unstable in $|D(T)|" \wedge (T$
superstable). 

Here the infinitary logic $\bbL_{\lambda^+,\kappa}$ is central.

A major point is to deal abstractly with what is essentially the
Boolean algebra of formulas over the empty set, $\bbB_T$ (so modulo
$T$ of course).  We introduce in Definition \ref{b4} the logics
$\bbL_{\lambda^+,\kappa}[\bbB]$ where $\bbB = \bbB_T$, the members of 
the Boolean algebra
(i.e. formulas from $\bbL(\tau_T)$) are coded by 
elements of the model and the union of these logics over the 
relevant $\bbB$'s is called $\bbL^{\ba}_{\lambda^+,\kappa}$, moreover
$\bbL^{\ba}_{\lambda,\kappa}$ is equivalent to 
$\bbL_{\lambda,\kappa}[\bbB^{\fr}_\lambda]$, see \ref{a4}(5).  
Then in Observation \ref{b5}(4) we note that:

\[
H(\bbL_{\lambda^+,\kappa}) \le H(\bbL_{\lambda^+,\kappa}[\bbB]) \le
H(\bbL^{\ba}_{\lambda^+,\kappa}) \le H(\bbL_{(2^\lambda)^+,\kappa}).
\]

\mn
The main result shows that there is an exact equivalence
between classes of the form $\bold N_{\lambda,T}$ and classes of the
form $\Mod_\psi,\psi \in \bbL_{\lambda^+,\kappa}[\bbB]$ for $\bbB$ 
the Boolean Algebra formulas over the emptyset in $T$.

We thank John Baldwin, Daniel Palacin and two referees for helpful comments.
\end{question}
\bigskip

\subsection {Preliminaries}\
\bigskip

Here for a first order complete $T$ we define the relevant parameters;
$\kappa(T),\bbB_T$ and quote characterization of the existence of
saturated models.

\begin{notation}
\label{a2}
1) $\tau$ will denote a vocabulary $\tau_M = \tau(M)$ is the vocabulary
of a model $M,|M|$ is the universe of $M$ and $\|M\|$ its cardinality;
$\bbL(\tau)$ is the first order logic for this vocabulary, i.e. the
set of first order formulas in $\tau$.

\noindent
1A) $T$ denotes a first order theory in $\bbL_{\tau(T)},\tau_T =
\tau(T)$ the vocabulary of $T$ and $T$ is complete and stable if not said
otherwise (but $T_1$ is neither necessarily complete nor necessarily stable).

\noindent
2) $\bar x_{[u]} = \langle x_i:i \in u\rangle$, similarly $\bar
   y_{[u]}$; e.g. $\bar x_{[\alpha]} = \langle x_i:i < \alpha\rangle$.

\noindent
3) $\bbL_{\lambda,\kappa}$ for $\lambda \ge \kappa$ is the logic where
the language $\bbL_{\lambda,\kappa}(\tau)$ is the following set of
formulas; it is the closure of the set of atomic formulas under negation,
conjunction of the form $\bigwedge\limits_{\alpha < \gamma}
\varphi_\alpha,\gamma < \lambda$ and quantification $(\exists 
\bar x_{[u]})\varphi$ where $u \in [\kappa]^{< \kappa}$ (really just
$(\exists \bar x_{[\varepsilon]})\varphi$ for $\varepsilon < \kappa$ suffice),
\but \, every formula has $< \kappa$ free variables.

\noindent
4) Let $\bbB$ denote a Boolean Algebra and $\uf(\bbB)$ the set of
   ultra-filters of $\bbB$.

\noindent
5) Let $\bold t$ denote an object as in Definition \ref{b2} below.

\noindent
6) For a theory $T$ let $\Mod_T$ be the class of models of $T$.
\end{notation}

\noindent
Recall
\begin{definition}  
\label{b13}
Let $T$ be a first order complete stable theory.

\noindent
0) For a model $M$ of $T$ and $A \subseteq M$ let $\bfS^n(A,M)$ be the
set of complete $n$-types over $A$ in $M$, equivalently $\{\tp(\bar
a,A,N):M \prec N$ and $\bar a \in {}^n N\}$ recalling that for $\bar a
\in {}^n M$ and $A \subseteq M$ we let $\tp(\omega,A,M) =
\{\varphi(\bar x_{[n]},\bar b + \varphi(\bar x,\bar y) \in
\bbL(\tau_M)$ and $\bar b \in {}^{\ell g(\bar y)}M$ and $M \models
\varphi[\bar a,\bar b]\}$; if $n=1$ then we may omit
$n$ and $\bfS^n(M) = \bfS^n(|m|,M)$ where $|M|$ is the universe of
$M$.

Recall:
\mn
\begin{enumerate}
\item[(a)]  $T$ is stable in $\lambda$ or $\lambda$-stable when for
  every model $M$ of $T$ and $A \subseteq M$ of cardinality $\le
  \lambda$ the set $\bfS(A,M)$ has cardinality $\le \lambda$
\sn
\item[(b)]  $T$ is superstable \Iff \, $T$ is $\lambda$-stable for
 every $\lambda$ large enough.
\end{enumerate}
\mn
1) $\kappa(T)$ is the minimal $\kappa$ such that: if $A \subseteq M_*
\in \Mod_T$ and $p \in \bold S(A,M)$ then there is $B \subseteq A$
of cardinality $< \kappa$ such that $p$ does not fork over $B$,
see \cite[Ch.III]{Sh:c}.

\noindent
2) Let $\kappa_r(T) = \min\{\kappa:\kappa$ regular $\ge \kappa(T)\}$
 so $\kappa_r(T)$ is the minimal regular
$\kappa$ such that $T$ is stable in $\lambda$ whenever $\lambda =
\lambda^{< \kappa} + 2^{|T|}$, see \cite[Ch.III]{Sh:c}.

\noindent
3) Let $\lambda(T)$ be the minimal $\lambda$ such that $T$ is stable
in $\lambda$, that is $[M \models T,\|M\| \le |T| + 
\aleph_0 \Rightarrow |\bold S(M)| \le \lambda]$,
 see \cite[Ch.III,\S5,\S6]{Sh:c}.

\noindent
4) $D_m(T) = \{\tp(\bar a,0,M):\bar a \in {}^m M \text{ and } M
\models T\}$ and $D(T) = \bigcup\limits_{m} D_m(T)$.

\noindent
5) Let $\EQ_T = \{\varphi(\bar x_{[n]},\bar y_{[n]}):n <
\omega,\varphi \in \bbL(\tau_T)$ and for every model $M$ of $T,\{(\bar
a,\bar b):\bar a,\bar b \in M$ and $M \models \varphi[\bar a,\bar
b]\}$ is an equivalence relation on ${}^n M$ with finitely many
equivalent classes$\}$.

\noindent
6) $M$ is $\aleph_\varepsilon$-saturated \when \, for every triple $(b,A,N)$
satisfying $A \subseteq M \prec N,b \in N,A$ finite, some $b' \in M$
realizes the type
$\{\varphi(x,b;\bar a):\bar a \subseteq A,\varphi(x,y,\bar a)$ is an
equivalence relation with finitely many equivalence classes in $M$,
this type is called $\stp(b,A,N)\}$, see \cite[Ch.III]{Sh:c}.
\end{definition}

\begin{remark}  
\label{b14}
By \cite[Ch.III,\S5,\S6]{Sh:c} we have that $\lambda(T) = 
|D(T)|^{<\kappa(T)}$ except when $|D(T)| < 2^{\aleph_0}$, if $T$ is
superstable and unstable in $|T|$, \then \,
$|D(T)| < 2^{\aleph_0} = \lambda(T)$ and $\lambda(T) = 
|D(T)|^{< \kappa(T)}$, see \ref{a12}.
\end{remark}

\noindent
The point is that by \cite[Ch.III]{Sh:c}:
\begin{fact}
\label{b18}
Let $T$ be a complete first order stable theory and let $\lambda \ge
\aleph_1 + |T|$ be an infinite cardinal.  \Then \, $T$ has a saturated
model of cardinality $\lambda$ if and only if $T$ is $\lambda$-stable,
if and only if $\lambda = \lambda^{< \kappa(T)} + \lambda(T)$.
\end{fact}

\noindent
Note that
\begin{observation}
\label{z6}
For every Boolean Algebra $\bbB_1$ of cardinality $\le
   \lambda$ and $\kappa \le \lambda^+$ there is a Boolean Algebra
$\bbB_2$ of cardinality $\lambda$ such that $|\uf(\bbB_2)| =
\Sigma\{|\uf(\bbB_1)|^\theta:\theta < \kappa\}$.
\end{observation}

\begin{PROOF}{\ref{z6}}
If $|\bbB_1| = \lambda,\kappa = \theta^+,\theta \le \lambda$ we
 define the Boolean Algebra $\bbB_2$ as the free product of 
$\theta$ copies of $\bbB_1$.

If $\kappa$ is a limit cardinal $\le \lambda,|\bbB_1|=\lambda$ let
   $\bbB_{2,\theta}$ be as above for $\theta < \kappa$ and $\bbB_2$
   the disjoint sum of $\langle \bbB_{2,\theta}:\theta <
   \kappa\rangle$ so essentially except one ultrafilter, all
   ultrafilters on $\bbB_2$ are ultrafilters on some $\bbB_{2,\theta}$
   so $\uf(\bbB_2) = 1 + \sum\limits_{\theta < \kappa} \uf(\bbB_{2,\theta})$.
\end{PROOF}

\begin{definition}  
\label{a4}
1) For a model $M$ and formula $\varphi(\bar x,\bar y) \in
 \bbL(\tau_M)$ and $\bar a \in {}^{\ell g(\bar y)}M$ let
 $\varphi(M,\bar a) = \{\bar b \in{}^{\ell g(\bar x)}M:M \models
 \varphi[\bar b,\bar a]\}$.

\noindent
2) For a model $M,\bbB_{M,m}$ is the Boolean Algebra of subsets of
   ${}^m M$ consisting of the sets $\{\varphi(M):\varphi =
   \varphi(\bar x_{[m]})\}$.

\noindent
2A) $\bbB_{T,m}$ for $T = \Th(M)$ is
the Boolean Algebra of the formulas $\varphi(\bar x_{[m]}) \in
\bbL(\tau_T)$ modulo equivalence
over $T$, so $\varphi_1(\bar x_{[m]}) \le \varphi_2(\bar x_{[m]})$
\Iff \, $T \vdash ``\varphi_1(\bar x_{[m]}) \rightarrow \varphi_2(\bar
x_{[m]})"$, so the elements are actually $\varphi(\bar x_{[m]})/\equiv_T$.

\noindent
3) Let $\bar{\bbB}_M = \langle \bbB_{M,m}:m < \omega\rangle$; abusing
notation let $\uf(\bar{\bbB}_M) =
\bigcup\limits_{m} \uf(\bbB_{M,m})$.  Similarly with $T$ instead of
$M$, also below.

\noindent
3A) Let $\bbB_M$ be the direct sum of $\langle \bbB_{M,m}:m <
\omega\rangle$ so $\langle 1_{\bbB_{M,m}}:m < \omega\rangle$ be a
maximal antichain of $\bbB_M,\bbB_M \rest \{x \in \bbB_M:x \le
1_{\bbB_{M,m}}\} = \bbB_{M,m}$ and $\cup \{\bbB_{M,m}:m < \omega\}$
generates $\bbB_M$.  
Let $\trufil(\bbB_M) =$ the
ultrafilter of $\bbB_M$ disjoint to $\{1_{\bbB_{M,n}}:n < \omega\}$ and let
$\uf^-(\bbB_M) = \uf(\bbB_M) \backslash \{\trufil(\bbB_M)\}$,
($\trufil$ stands for trivial ultra-filter).

\noindent
4) Let $\lambda'(M)$ be the cardinality of $\uf(\bbB_M)$ and
$\lambda'(T) = \lambda'(M)$ when $M \models T$.

\noindent
5) Let $\bbB^{\fr}_\lambda$ be the Boolean algebra generated freely by
$\{\bold a_\alpha:\alpha < \lambda\}$ so $\uf(\bbB^{\fr}_\lambda)$ has
cardinality $2^\lambda$.
\end{definition}

\begin{remark}
\label{a5}
We may be interested in the Boolean Algebra of formulas which are
almost over $\emptyset$, i.e. $\varphi(\bar x_m,\bar a),\bar a \in
{}^{\ell g(\bar y)}M$ where $\varphi(\bar x_m,\bar y) \in \bbL(\tau_T)$
satisfies: $\varphi(\bar x_m,\bar y)$ such that for some 
$\vartheta(\bar x_m,\bar y_m) \in \EQ^m_M$, see \ref{b13}(5), we have $M
\models (\forall \bar z)(\forall \bar x_m,\bar y_m)[\vartheta(\bar
  x_m,\bar y_m) \rightarrow (\varphi(\bar x_m,\bar z) \equiv
\varphi_n(\bar y_m,\bar z)]$.

But this is not necessary here.
\end{remark}

\begin{observation}
\label{a6}
1) $\bbB_{M,m}$ essentially depend just on $\Th(M)$, i.e. if $T =
 \Th(M)$ then $\bbB_{M,m}$ is isomorphic to $\bbB_{T,m}$ where an
 isomorphism $\bfj$ is defined as follows:
$\varphi(\bar x_{[m]}) + \bbL(\tau_T) \Rightarrow \bold j(\varphi(M))
 = \varphi(\bar x_{[m]})/\equiv_T$, so $\lambda'(T)$ is well defined.

\noindent
2) Similarly for other notions from Definition \ref{a4}.

\noindent
3) $\uf^-(\bbB_M),\uf(\bbB_M)$ has the same cardinality, in fact,
there is a natural one-to-one mappping $\pi$ from $\uf(\bar{\bbB}_M)$
onto $\uf^-(\bbB_M)$ such that $D \in \uf(\bbB_{M,m}) \Rightarrow
\pi(D) = \{a \in \bbB_{M,m}:a \cap 1_{\bbB_{M,m}} \in D\}$.
\end{observation}

\noindent
Recall by Lemma \cite[Ch.III,3.10]{Sh:c}:
\begin{fact}
\label{a9}
Let $T$ be a stable (first order complete) theory, $\kappa = \kappa(T)$ 
and $M$ is an uncountable model of $T$.  Then $M$ is saturated \Iff \,
\medskip

\noindent
\underline{Case 1}: $\kappa > \aleph_0$
\mn 
\begin{enumerate}
\item[$(a)$]  if $\bold I \subseteq M$ is an infinite indiscernible
  set \then \, there is an indiscernible set $\bold J \subseteq M$
  extending $\bold I$ of cardinality $\|M\|$
\sn
\item[$(b)$]  $M$ is $\kappa$-saturated.
\end{enumerate}
\medskip

\noindent
\underline{Case 2}: $\kappa = \aleph_0$
\mn 
\begin{enumerate}
\item[$(a)'$]  if $A \subseteq M$ is finite and $a \in M \backslash
 \acl(A)$ \then \, there is an indiscernible set $\bold J$ over $A$ in
 $M$ based on $A$ such that $a \in \bold J$ and $\bold J$ is 
of cardinality $\|M\|$
\sn
\item[$(b)'$]  $M$ is $\aleph_\varepsilon$-saturated, see \cite{Sh:c}
  or Definition \ref{b13}(6).
\end{enumerate}
\end{fact}

\begin{fact}
\label{a12}
Assume $T$ is a stable (first order complete) theory.

\noindent
1) If $\kappa(T) > \aleph_0$ then $\lambda(T) =
   |D(T)|^{<\kappa_r(T)}$.

\noindent
2) If $\kappa(T) = \aleph_0$ \then \, $\lambda(T)$ is $|D(T)|$ or
$\lambda(T) = 2^{\aleph_0} + |D(T)|$ and 
\mn
\begin{enumerate}
\item[$(\st)_T$]  for some finite $A \subseteq M,M\in
   \Mod_T$, the set $\{\stp(a,A):a \in M\}$ has cardinality continuum.
\end{enumerate}
\end{fact}

\begin{definition}
\label{a16}
1) For a cardinal $\theta$ let $T^{\eq}_\theta$ be the model
completion of $T^{\eq,0}_\theta$, see below.

\noindent
2) Let $\tau^{\eq}_\theta = \{E_i:i < \theta\},E_i$ a two-place predicate.

\noindent
3) Let $T^{\eq}_\theta$ be the universal theory included in
$\bbL(\tau^{\eq}_\theta)$ such that: for a $\tau^{\eq}_\theta$-model
$M,M \models T^{\eq}_\theta$ \Iff \, $E^M_i$ is an equivalence
relation and $E^M_j$ refines $E^M_i$ for $i < j < \theta$.
\end{definition}

\begin{claim}
\label{a19}
(Basic properties of non-forking)

\noindent
1) $M_\delta = \bigcup\limits_{i < \delta} M_i$ is $\lambda$-saturated
\when \,:
\mn
\begin{enumerate}
\item[(a)]  $\langle M_i:i < \delta\rangle$ is a $<$-increasing
  sequence of models of $T$
\sn
\item[(b)]  $T$ is stable and $\kappa(T) \le \cf(\delta)$
\sn
\item[(c)]  each $M_i$ is $\lambda$-saturated.
\end{enumerate}
\mn
2) If $T$ is superstable, $\lambda(T) > |D(T)|$ - FILL.
\end{claim}

\begin{PROOF}{\ref{a19}}
1) See \cite[Ch.III]{Sh:c}.

\noindent
2) See \cite[Ch.III,5.9,5.10,5.11]{Sh:c}.
\end{PROOF}
\newpage

\section {The frame}

First, we define here $\bold N_{\lambda,T}$, the set of triples $\bold t$ from
the abstract when we fix $T,\lambda$ and for $\bold t \in \bold N_{\lambda,T}$
we define the class of models $\Mod_{\bold t}$ (in \ref{b2},\ref{b2g})
and give easy properties (in \ref{b3}, \ref{b3f}).  Second, we deal with
the logics $\bbL_{\lambda,\kappa}[\bbB]$ via which we shall
characterize the Hanf number of $\bold N_{\lambda,T}$ and look at the
relations among such logics (see \ref{b4}, \ref{b10}, \ref{b12}).
Third, we deal
with representations, e.g. how $\psi \in \bbL_{\lambda^+,\kappa}$ can
be translated to models of first order $T$, with extra demands (see
\ref{b14} - \ref{b24}).  Lastly, we look at order between the $\bbB$'s.

\begin{definition}  
\label{b2}
1) For $T$ complete first order stable theory and $\lambda \ge |T|$ let 
$\bold N_{\lambda,T}$ be the class of triples $\bold t = (T,T_1,p) =
(T_{\bold t},T_{1,\bold t},p_{\bold t})$ such that:
\mn
\begin{enumerate}
\item[$(a)$]  $T_{\bold t} = T$
\sn
\item[$(b)$]  $T_1 \supseteq T$ is a first order theory and
  $|\tau(T_1)| \le \lambda$
\sn
\item[$(c)$]  $p(x)$ is an $\bbL(\tau_{T_1})$-type, not necessarily complete.
\end{enumerate}
\mn
1A) For $\bold t$ as above we say $M_1 \models \bold t$ or $M_1 \in
\Mod_{\bold t}$ or $M_1$ is a model of $\bold t$ \when \,:
\mn
\begin{enumerate}
\item[$(a)$]  $M_1 \models T_{1,\bold t}$ and $M_1$ a $\tau_{T_1}$-model 
\sn
\item[$(b)$]  $M_1$ omits the type $p_{\bold t}(x)$
\sn
\item[$(c)$]  $M_1 \rest \tau_T$ is saturated.
\end{enumerate}
\mn
1B) Omitting $T$ means: for some $T$.

\noindent
2) Let $\spec_{\bold t} = \{\|M\|:M \models \bold t\}$ 
for $\bold t \in \bold N_{\lambda,T}$.

\noindent
3) The Hanff number 
$H(\bold N_{\lambda,T})$ is the minimal $\mu$ such that: if $\bold
   t \in \bold N_{\lambda,T}$ and $\bold t$ has a model of cardinality
   $\ge \mu$ \then \, $\bold t$ has models of arbitrarily large
   cardinality; see \ref{b4}(3).

\noindent
3A) Equivalently, $H(\bold N_{\lambda,T}) = \sup\{H(\bold t):H(\bold t) <
\infty,\bold t \in \bold N_{\lambda,T}\}$ where $H(\bold t) =
\sup\{\|M\|^+:M \in \Mod_{\bold t}\}$.

\noindent
4) $\lambda(\bold t) := \lambda(T_{\bold t}) + |T_{1,\bold t}|$
recalling \ref{b13}(3).
\end{definition}

\begin{convention}
\label{b2g}
Below $\bold t,T,T_1,p,\lambda$ are as in
Definition \ref{b2} if not said otherwise and then $\kappa =
\kappa_r(T)$ is as in \ref{b13}.
\end{convention}

\begin{claim}
\label{b3}
1) If $M \in \Mod_{\bold t}$ has cardinality $\mu$ \then \, $\mu = \mu^{<
\kappa(T)} + |\lambda(T)|$, i.e. $\mu \in \spec_{\bold t}
   \Rightarrow \mu = \mu^{< \kappa(T)} + \lambda(T)$.

\noindent
2) If $M \in \Mod_{\bold t}$ and $\lambda(\bold t) \le \mu = \mu^{<
  \kappa(T)} < \|M\|$ recalling \ref{b2}(4) 
and $A \subseteq M$ is of cardinality $\mu$
 \then \, for some $N$ we have:
\mn
\begin{enumerate}
\item[$(a)$]  $N \in \Mod_{\bold t}$
\sn
\item[$(b)$]  $A \subseteq N \prec M$
\sn
\item[$(c)$]  $N$ has cardinality $\mu$.
\end{enumerate}
\end{claim}

\begin{PROOF}{\ref{b3}}
1) By \ref{b18}.

\noindent
2) Note that also $\mu = \mu^{< \kappa_r(T)}$ by cardinal arithmetic and
hence $\kappa_r(T) \le \mu$; we choose $M_i$ by induction on 
$i < \kappa_r(T)$ such that:
\mn
\begin{enumerate}
\item[$(a)$]  if $i$ is even then $M_i \prec M$ and $\|M_i\| = \mu$
\sn
\item[$(b)$]  if $i$ is odd then $M_i \rest \tau(T_{\bold t}) 
\prec M \rest \tau(T_{\bold t}),
\|M_i\| = \mu$ and $M_i$ is saturated
\sn
\item[$(c)$]  if $j<i$ then $A \cup |M_j| \subseteq |M_i|$.
\end{enumerate}
\mn
There is no problem to carry the induction and then $M' = \cup\{M_{2i}:i <
\kappa_r(T)\} = \cup\{M_{2i+1}:i < \kappa_i(T)\}$ 
is as required: $M' \prec M$ by (a)+(c) and
Tarski-Vaught, $\|M'\| = \mu$ since $\mu^{< \kappa_T(T)} = \mu$
and $M' \rest \tau(T)$ is saturated by (b) + (c) and \ref{a19}(2).
\end{PROOF}

\begin{conclusion}
\label{b3f}
For understanding the Hanf number of $\bold t$, it is enough to 
consider cardinals $\mu = \mu^{< \kappa(T)} \ge \lambda(\bold t)$.
\end{conclusion}

\noindent
Now we turn to the logics of the form $\bbL_{\lambda^+,\kappa}[\bbB]$;
first we define them.
\begin{definition}
\label{b4}
1) Assume
\mn
\begin{enumerate}
\item[$(a)$]   $\lambda \ge \kappa = \cf(\kappa)$
\sn
\item[$(b)$]  $\bbB$ is a Boolean Algebra of cardinality $\le \lambda$ and
recall $\uf(\bbB)$ is the set of ultrafilters on $\bbB$.
\end{enumerate}
\mn
\Then
\mn
\begin{enumerate}
\item[$(\alpha)$]  Let $\voc_\lambda[\bbB]$ be the class of
vocabularies $\tau$ of cardinality $\le \lambda$ such that 
$c_b \in \tau$, an individual constant for each $b \in \bbB$,
and $P,Q \in \tau$ unary predicates and $R \in \tau$ binary predicate
and $\tau$ may have additional signs.
\sn
\item[$(\beta)$]   For $\tau \in \voc_\lambda[\bbB]$ let 
$\bbL_{\lambda^+,\kappa}[\bbB](\tau)$ be the 
set of sentences $\psi \in \bbL_{\lambda^+,\kappa}(\tau)$ but we
stipulate that from $\psi$ we can reconstruct the triple
$(\lambda^+,\kappa,\bbB)$ hence $\bbL_{\lambda^+,\kappa}[\bbB]$.
\end{enumerate}
\mn
[Note that $\psi$ has $\le \lambda$ sub-formulas]:
\mn
\begin{enumerate}
\item[$(\gamma)$]   omitting $\tau$ means $\tau = \tau_\psi$ is 
the minimal $\tau \in \voc_\lambda[\bbB]$ such that
$\psi \in \bbL_{\lambda^+,\kappa}[\bbB](\tau)$.
\end{enumerate}
\mn
2) For $\tau \in \voc_\lambda[\bbB]$ and 
$\psi \in \bbL_{\lambda^+,\kappa}[\bbB](\tau)$ let
$\Mod^1_\psi[\bbB]$ be the class of models $M$ of $\psi$ (which are
$\tau_\psi$-models if not said otherwise) such that (note: clauses
(a)-(e) can be expressed in $\bbL_{\lambda^+,\aleph_0}$, but when
$|\uf(\bbB)| > \lambda$ not so clause (f)):
\mn
\begin{enumerate}
\item[$(a)$]  $P^M = \{c^M_b:b \in \bbB\}$
\sn
\item[$(b)$]  $\langle c^M_b:b \in \bbB\rangle$ are pairwise distinct
\sn
\item[$(c)$]  $R \subseteq P^M \times Q^M$
\sn
\item[$(d)$]  for every $a \in Q^M$ the set
  $\uf^M(a) := \{b \in \bbB:M \models c_b Ra\}$ belongs to $\uf(\bbB)$
\sn
\item[$(e)$]  if $a_2 \ne a_2$ are from $Q^M$ then $\uf^M(a_1) \ne
  \uf^M(a_2)$
\sn
\item[$(f)$]  for every $u \in \uf(\bbB)$ there is $a \in Q^M$ 
such that $M \models \bigwedge\limits_{i < \lambda}
(c_b R a)^{\iif(b \in u)}$, (by clause (e) the element $a$ is unique).
\end{enumerate}
\mn
3) Let $\Mod^2_\psi[\bbB]$ be the class of $M \in \Mod^1_\psi[\bbB]$ such that:
\mn
\begin{enumerate}
\item[$(f)$]  $\|M\| = \|M\|^{< \kappa}$ and (follows) $\|M\| \ge |\uf(\bbB)|$.
\end{enumerate}
\mn
4) For $\iota = 1,2$ and $\psi \in \bbL_{\lambda^+,\kappa}[\bbB]$ let
$\spec^\iota_\psi[\bbB] = \{\|M\|:M \in \Mod^\iota_\psi[\bbB]\}$.

\noindent
4A) Writing $\Mod^\iota_\psi,\spec^\iota_\psi$ we mean $\iota \in
\{1,2\}$ and may omit $\iota$ when $\iota=2$ (because this is the main
case for us), see \ref{b5}(1) below 
and $\bbB$ can be reconstructed from $\psi$.

\noindent
5) Let $H(\bbL_{\lambda^+,\kappa}[\bbB])$ be the first $\mu$ such 
that: if $\psi \in \bbL_{\lambda^+,\kappa}[\bbB]$ and there is $M \in 
\Mod_\psi[\bbB]$ of cardinality $\ge \mu$ \then \, $\{\|M\|:M \in 
\Mod_\psi[\bbB]\}$ is an unbounded class of cardinals.

\noindent
6) Let $\bbL^{\ba}_{\lambda^+,\kappa}$ be
$\cup\{\bbL_{\lambda^+,\kappa}[\bbB]: \bbB$ a Boolean\footnote{So
  every sentence $\psi \in \bbL^{\ba}_{\lambda^+,\kappa}$ fixes a
  Boolean Algebra $\bbB$ as above and a vocabulary of cardinality
 $\le \lambda$ from $\voc_\lambda[\bbB]$ as described.} Algebra 
of cardinality $\le \lambda\}$ so
every sentence of $\bbL^{\ba}_{\lambda^+,\kappa}(\tau)$ is a sentence in
$\bbL_{\lambda^+,\kappa}[\bbB](\tau)$ for some $\bbB$ as above; 
so we may stipulate that the set of
elements of $\bbB$ is a cardinal $\le \lambda$ and $c_i\in \tau$ for
$i < \lambda$.

\noindent
7) We define $H(\bbL^{\ba}_{\lambda^+,\kappa})$ similarly; yes, this
   is just $\sup\{H(\bbL_{\lambda^+,\kappa}[\bbB]):\bbB$ as above$\}$.
\end{definition}

\noindent
Having defined the sets $(\bbL_{\lambda^+,\kappa}[\bbB])(\tau)$ of
sentences and the relevant classes of models $\Mod^\iota_\psi[\bbB]$ and
spectrums $\spec^\iota_\psi[\bbB]$ and Hanf numbers we should now try to
understand the order between them.
\begin{claim}
\label{b6}
1) Recalling $\bbB^{\fr}_\lambda$ is the Boolean Algebra generated
freely by $\lambda$ generators:
\mn
\begin{enumerate}
\item[$(a)$]  for every Boolean algebra $\bbB_1$ of cardinality
  $\lambda$ or just $\le \lambda$ 
and $\psi_1 \in \bbL_{\lambda^+,\kappa}[\bbB_1]$ there is
  $\psi \in \bbL_{\lambda^+,\kappa}[\bbB^{\fr}_\lambda]$ such that
$\spec^{\iota}_{\psi_1} \backslash 2^\lambda = 
\spec^{\iota}_\psi \backslash 2^\lambda$ for $\iota= 1,2$
\sn
\item[$(b)$]  $H(\bbL_{\lambda^+,\kappa}[\bbB_1]) \le
H(\bbL_{\lambda^+,\kappa}[\bbB^{\fr}_\lambda])$ for $\bbB_1$ as above.
\end{enumerate}
\mn
2) If $\bbB_1,\bbB_2$ are Boolean algebras of cardinality $\le
\lambda$ and $\bbB_1$ is a homomorphic image of $\bbB_2$, \then \,:
\mn
\begin{enumerate}
\item[$(a)$]  for every $\psi_1 \in \bbL_{\lambda^+,\kappa}[\bbB_1]$ there is
$\psi_2 \in \bbL_{\lambda^+,\kappa}[\bbB_2]$ such that
$\spec^{\iota}_{\psi_1}[\bbB_1] = \spec^{\iota}_{\psi_2}[\bbB_1]$ 
for $\iota= 1,2$
\sn
\item[$(b)$]  $H(\bbL_{\lambda^+,\kappa}[\bbB_1]) \le
H(\bbL_{\lambda^+,\kappa}[\bbB_2]$.
\end{enumerate}
\mn
3) For every $\psi_1 \in \bbL_{\lambda^+,\kappa}[\bbB]$ there
 are $\psi_2,\psi'_2,\psi''_2 \in \bbL_{\lambda^+,\kappa}[\bbB]$ such that:
\mn
\begin{enumerate}
\item[$(a)$]   $\spec^1_{\psi_2}[\bbB] = \{\mu:\mu =\mu^{< \kappa} \in
\spec^1_{\psi_1}[\bbB]\} = \spec^2_{\psi_1}[\bbB]$ 
and\footnote{Recall that if $\mu > 2^{<
 \kappa}$ then $(\mu^{< \kappa})^{< \kappa} = \mu$, see \cite{Sh:g}.}
\sn
\item[$(b)$]   $\spec^1_{\psi'_2}[\bbB] =
 \{\mu^{< \kappa}:\mu \in \spec^1_{\psi_1}[\bbB]\}$ and
\sn
\item[$(c)$]   $\spec^1_{\psi''_2}[\bbB] = \{\mu:\mu \ge \lambda$ and
  $\mu \in \spec^1_{\psi_1}[\bbB]\}$.
\end{enumerate} 
\end{claim}

\begin{PROOF}{\ref{b6}}
1) Let $h$ be a homomorphism from $\bbB^{\fr}_\lambda$ onto $\bbB_1$,
exists as $\bbB_1$ is a Boolean algebra of cardinality $\le \lambda$.
Now apply part (2).

\noindent
2) Let $I := \Ker(h) := \{a \in \bbB^{\fr}_\lambda:h(a)=0\}$ and let
$h_1:\bbB_2 \rightarrow \bbB_2$ be such that $a \in \bbB_2 \Rightarrow
h(h_2(a)) = a$.  Let $\bbB'_1$ be the Boolean Algebra with set of
elements $\Rang(h_1)$ such that $h_2$ is an isomorphism from $\bbB_1$
onto $\bbB'_1$.  Let $\psi'_1$ be like $\psi_1$ replacing $\bbB_1$ by
$\bbB'_1$ and the predicate $P$ by a predicate $P'$.  The rest should
be clear.

\noindent
3) Should be clear but we elaborate.
\medskip

\noindent
\underline{Clause (a)}:  Let $\tau_2 = \tau(\psi_1) \cup \{F_{i,j}:i <
j < \kappa\}$ with $F_{i,j} \notin \tau(\psi)$ be pairwise distinct
unary function.

Let $\psi_2 = \psi_1 \wedge \varphi_2$ where

\[
\varphi_2 = \bigwedge\limits_{0 <j<\kappa} (\forall
\ldots,x_i,\ldots)_{i<j} (\exists y)[\bigwedge\limits_{i<j} F_i(y) =
x_i].
\]

\mn
Now think
\medskip

\noindent
\underline{Clause (b)}:  Let $\tau'_2 = \tau(\psi_1) \cup \{F_{i,j}:i <
j < \kappa\} \cup \{P_j:j<\kappa\}$ with $F_{i,j}$ as above $P_j,P 
\notin \tau(\psi_1)$ be pairwise distinct unary predicate.

Let $\psi^p_1 \in \bbL_{\lambda^+,\kappa}[\bbB]$ be  such that for a
$(\tau(\psi_1) \cup \{P\})$-model $M,M \models \psi^p_1$ \Iff \, $(M
\rest P^M) \rest \tau(\psi_1)$ is a $\tau(\psi_1)$-model and is a
model of $T$.

Lastly, let $\psi_2 = \psi^p_1 \wedge \varphi'_2$ where $\varphi'_2$
is the conjunction of:
\mn
\begin{itemize}
\item  $M \models \varphi^0_2$ \Iff \, $\langle P^M\rangle \char 94
  \langle P^M_j:j<\kappa\rangle$ is a partition of $M$
\sn
\item  $\varphi^1_{2,i,j} = (\forall x)(P(F_{i,j}(x))$ for $i < j <
  \kappa$
\sn
\item  $\varphi^2_{2,j} = (\forall x,y)[x \ne y \wedge P_j(x)
  \rightarrow \bigvee\limits_{i<j} F_{i,j}(x) \ne F_{i,j}(y)]$
\sn
\item  $\varphi^3_{2,j} = (\forall \ldots,x_i,\ldots)_{i<j}
(\bigwedge\limits_{i<j} P(x_i) \rightarrow (\exists y)(P_j(y) \wedge
\bigwedge\limits_{i<j} F_{i,j}(y) = x_i)$.
\end{itemize}
\mn
Now check.
\medskip

\noindent
\underline{Clause (c)}:  

Even easier.
\end{PROOF}

\begin{observation}  
\label{b5}
Let $\bbB$ be a Boolean Algebra of cardinal $\le \lambda$ and $\kappa
\le \lambda^+$.

\noindent
1) In the Definition \ref{b4}(5) of $H(\bbL_{\lambda^+,\kappa}[\bbB])$
it does not matter if we use $\Mod^1_\psi[\bbB]$ or
$\Mod^2_\psi[\bbB]$.

\noindent

\noindent
2) For every $\mu < H(\bbL_{\lambda^+,\kappa}[\bbB])$ we have $2^\mu < 
H(\bbL_{\lambda^+,\kappa}[\bbB])$ hence
$H(\bbL_{\lambda^+,\kappa}[\bbB])$ is a strong limit cardinal of
cofinality $> \lambda$.

\noindent
3) $H(\bbL^{\ba}_{\lambda^+,\kappa}) <
H(\bbL_{(2^\lambda)^+,\kappa})$.

\noindent
4) We have $H(\bbL_{\lambda^+,\kappa}) \le 
H(\bbL_{\lambda^+,\kappa}[\bbB]) \le
\bbL_{\lambda^+,\kappa}[\bbB^{\fr}_\lambda] = H(\bbL^{\ba}_{\lambda^+,\kappa})
< H(\bbL_{(2^\lambda)^+,\kappa})$.

\noindent
5) If $\bbB^{\fr}_\lambda$ is the free Boolean Algebra of cardinality
$\lambda$ from \ref{a4}(5) and 
$\kappa = \aleph_0$ \then \, $H(\bbL_{\lambda^+,\kappa}) <
\beth_{(2^\lambda)^+} <
H(\bbL_{\lambda^+,\kappa}[\bbB^{\fr}_\lambda])$.  Also for any $\kappa
\ge \aleph_0$ we have $H(\bbL^{\ba}_{\lambda^+,\kappa}) <
H(\bbL_{(2^\lambda + \aleph_0})$.

\noindent
6) If $\psi \in \bbL_{\lambda^+,\kappa}[\bbB]$ and
 $H(\bbL_{\lambda^+,\kappa}[\bbB]) \le \sup\{\|M\|:
M \in \Mod_\psi[\bbB]\}$ \then \, $\infty = \sup\{\|M\|:
M \in \Mod_\psi[\bbB]\}$ hence $\cf(H(\bbL_{\lambda^+,\kappa}[\bbB]))
\le 2^\lambda$).

\noindent
7) Like part (5) for $\psi \in \bbL^{\ba}_{\lambda^+,\kappa}$ and
$\Mod^{\ba}_\psi$.
\end{observation}

\begin{PROOF}{\ref{b5}}
1) First, as easily the Hanf number is $> 2^\lambda \ge |\uf(\bbB)|$,
we can ignore models of cardinality $< 2^\lambda$.
Second,
\mn
\begin{enumerate}
\item[$(*)_1$]   if $\psi_1 \in \bbL_{\lambda,\kappa}[\bbB](\tau)$ and
 $\sup(\spec^1_{\psi_1}) < \infty$ then 
$\sup(\spec^2_{\psi_1}) \le \sup(\spec^1_{\psi_1}) \le
(\sup(\spec^2_{\psi_1}))^{< \kappa} < \infty$.
\end{enumerate}
\mn
[Why? the first inequality
because $\spec^1_\psi \supseteq \spec^2_\psi$; the second inequality
by \ref{b3}(2).]

We can conclude that the Hanf number of the logic
$\bbL_{\lambda^+,\kappa}[\bbB]$ using $\Mod^1_\psi$ is smaller or
equal to the Hanf number of the logic $\bbL_{\lambda^+,\kappa}[\bbB]$
using $\Mod^2_\psi$.  Alternatively, if $\psi_1 \in
\bbL_{\lambda^+,\kappa}[\bbB]$ then by \ref{b6}(3)(b) there is
$\psi'_2 \in \bbL_{\lambda^+,\kappa}[\bbB]$ such that
$\sup(\spec^1_{\psi_1}) < \infty \Rightarrow \sup(\spec^1_{\psi_1}) \le
\sup(\spec^2_{\psi'_2}) < \infty$, hence the Hanf number using
$\spec^1_\psi$'s is $\le$ the Hanf number using $\spec^2_\psi$'s.
Moreover, above we get $\sup(\spec^1_{\psi_1}) \le
\sup(\spec^2_{\psi'_2}) = \sup(\spec^1_{\psi'_2})$ as
$\spec^2_{\psi'_2} = \spec^1_{\psi'_2}$.

On the other hand, by clause (a) of \ref{b6}(3) if $\psi_1 \in
\bbL_{\lambda,\kappa}[\bbB]$ then there is $\psi_2 \in
\bbL_{\lambda,\kappa}[\bbB]$ such that $\spec^1_{\psi_2} =
\spec^2_{\psi_1}$ so $\sup(\spec^2_{\psi_1}) < \infty \Rightarrow \sup
\spec^2_{\psi_1} = \sup \spec^1_{\psi_2} < \infty$
so also the other inequality holds.

\noindent
2) For any $\psi_1 \in \bbL_{\lambda^+,\kappa}[\bbB]$ we can find $\psi_2 \in
\bbL_{\lambda^+,\kappa}[\bbB]$ such that $\tau_{\psi_1} \subseteq
   \tau_{\psi_2},P_*,R_* \in \tau_{\psi_2} \backslash \tau_{\psi_1}$ 
are unary, binary predicates respectively and:
\mn
\begin{enumerate}
\item[$(*)_1$]  $M_2 \in \Mod^\iota_{\psi_2}[\bbB]$ \Iff \,
\sn
\begin{itemize}
\item  $(M_2 \rest P^{M_2}_* \rest \tau_{\psi_1}) \in \Mod_{\psi_1}[\bbB]$
\sn
\item  $M_2 \models (\forall y,z)(\exists
  x)[P_*(x) \wedge (R(x,y) \equiv \neg R(x,z))]$ hence 
$|P^{M_2}_*| \le \|M_2\| \le 2^{|P_*(M_2)|}$.
\end{itemize}
\end{enumerate}
\mn
Clearly
\mn
\begin{enumerate}
\item[$(*)_2$]  for every $M_1 \in \Mod^1_{\psi_1}[\bbB]$ and $\mu
  = \mu^{< \kappa} \in [\|M_1\|,2^{\|M_1\|}]$ there is $M_2 
\in \Mod^1_{\psi_2}[\bbB]$ of cardinality $\mu$.
\end{enumerate}
\mn
Using $(*)_2$ this clearly suffices for the first statement.  The
second is easy, too.

\noindent
3)  Let $\bold K_{\lambda^+,\kappa}$ be the class of 
pairs $(\psi,\bbB)$ such that $\bbB$ is a Boolean Algebra of cardinality $\le
\lambda,\psi \in \bbL_{\lambda^+,\kappa}[\bbB]$.  For $(\psi,\bbB)
\in \bold K_{\lambda^+,\kappa}$ let $H(\psi,\bbB) = \cup\{\mu^+:\mu
\in \spec^2_\psi(\bbB)\}$.  Clearly up to
isomorphism (of vocabularies) $\bold K_{\lambda^+,\kappa}$ has
cardinality $\le 2^\lambda$ and hence $\bold C_{\lambda^+,\kappa} :=
\{H(\psi,\bbB):(\psi,\bbB) \in \bold K_{\lambda^+,\kappa}\}$ has
cardinality $\le 2^\lambda$.  So let $\langle
(\psi_i,\bbB_i):i < 2^\lambda\rangle$ be such that $(\psi_i,\bbB_i)$
is as above and $\bold C_{\lambda^+,\kappa} \backslash \{\infty\} =
\{\mu_i:i < 2^\lambda\}$ where $\mu_i = H(\psi_i,\bbB_i) = \cup\{\mu^+:\mu \in
\spec^1_{\psi_i}[\bbB_i]\}$ for $i < 2^\lambda$.  Now we can find $\psi \in
\bbL_{(2^\lambda)^+,\kappa}$ such that $M \models \psi$ iff
\mn
\begin{enumerate}
\item[$(*)$]   $<^M$ is a linear order of $|M|$ and for arbitrarily
large $a \in M$ there are $i < 2^\lambda$ and 
$N \in \Mod^2_{\psi_i}[\bbB_i]$ with universe $\{b:b <^M a\}$.
\end{enumerate}
\mn
Together with part (2), clearly $\infty > \sup(\spec_\psi) =
\max(\spec_\psi) = \cup\{\mu_i:i < 2^\lambda\}$ so we are done.

\noindent
4) For the first inequality ``$H(\bbL_{\lambda^+,\kappa}) \le
H(\bbL_{\lambda^+,\kappa}[\bbB])$", see the definitions of
$\bbL_{\lambda^+,\kappa}[\bbB]$.  For the second inequality,
``$H(\bbL_{\lambda^+,\kappa}[\bbB]) \le
H(\bbL_{\lambda^+,\kappa}[\bbB^{\fr}_\lambda])$", use \ref{b6}(1)(b).  For
the third inequality, ``$H(\bbL_{\lambda^+,\kappa}[\bbB^{\fr}_\lambda]) =
H(\bbL^{\ba}_{\lambda^+,\kappa})$", use the definition of the latter
and the second inequality.  For the fourth inequality,
``$H(\bbL^{\ba}_{\lambda^+,\kappa}) <
H(\bbL_{(2^\lambda)^+,\aleph_0})$", the inequality holds as every model
$M$ satisfying $\lambda \ge \|M\| + |\tau_M|$ can be characterized up
to isomorphism by some $\psi \in \bbL_{(2^\lambda)^+,\kappa}$.

\noindent
5) The first inequality ``$H(\bbL_{\lambda^+,\kappa}) <
\beth_{(2^\lambda)^+}$" holds, is well known see, e.g. by 
Theorem 5.4 and 5.5 of \cite[Ch.VII,\S5]{Sh:c}
recalling $\kappa = \aleph_0$.  The second inequality, 
``$\beth_{(2^\lambda)^+} <
H(\bbL_{\lambda^+,\kappa}[\bbB^{\fr}_\lambda])$", holds by the
equality in part (4) and part (3).

For the third inequality note that:
\mn
\begin{enumerate}
\item[$(*)$]  there is $\psi \in
\bbL_{\lambda^+,\kappa}[\bbB^{\fr}_\lambda]$ such that: $M \models
\psi$ \Iff \,:
\sn
\begin{enumerate}
\item[(a)]  $P^M,Q^M,R^M$ are as in Definition
\sn
\item[(b)]  $F^M_i (i < \lambda)$ are as in \ref{b6}(3)(a) for $Q^M$,
  i.e. $M \models (\forall \ldots,x_i,\ldots)_{i < \lambda} (\exists
  y)[\bigwedge\limits_{i < \lambda} Q(x_i) \rightarrow (\exists
  y)(Q(y) \wedge \bigwedge\limits_{i < \lambda} F_i(y) = x_i)]$
\sn
\item[(c)]  $<^M$ is a well ordering of $Q^M$.
\end{enumerate}
\end{enumerate}
\mn
6) As in the end of the proof of part (3) replacing $\psi_i$ by
$\psi$, that is, we can find $\psi_i \in
\bbL_{\lambda^+,\kappa}[\bbB]$ such that:
\mn
\begin{enumerate}
\item[$(*)$]  $M_1 \models \psi_1$ \Iff \, for some $< \in
  \tau(\psi_1),<^{M_1}$ is a linear order of $(M_1)$ such that for
  arbitrarily large $b \in M_1,M_1 \rest \{a:a <^{M_1}\} \rest
  \tau_\psi$ is a model of $\psi$.
\end{enumerate}
\mn
Clearly this suffice.

\noindent
7) So assume $\mu < H(\bbL^{\ba}_{\lambda^+,\kappa})$ hence by the
definition there is $\psi \in \bbL^{\ba}_{\lambda^+,\kappa}$ such
that $\{\|M\|:M \models \psi\}$ is bounded by as a member $\ge
  \mu$.  By the definition of $\bbL^{\ba}_{\lambda^+,\kappa}$ for some
  Boolean Algebras $\bbB$ of cardinality $\le \lambda$ we have $\psi
  \in L^{\ba}_{\lambda^+,\kappa}[\bbB]$ and now apply part (2).
\end{PROOF}

\noindent
The following \ref{b9}, \ref{b10}, \ref{b12} is another way to represent the
logic $\bbL^{\ba}_{\lambda^+,\kappa}$ equivalently the logic
$\bbL_{\lambda^+,\kappa}[\bbB^{\fr}_\lambda]$, hence eventually
to state the Hanf numbers.
\begin{definition}  
\label{b9}
1) Let $\bbL^*_{\lambda^+,\kappa}$ be defined like
$\bbL^{\ba}_{\lambda^+,\kappa}$, see \ref{b2}(3) replacing $\langle
   c_b:b \in \bbB\rangle$ by $\langle c_i:i < \lambda\rangle$ and
   $\uf(\bbB)$ by $\cP(\{c_i:i < \lambda\})$.

\noindent
2) For $\psi \in \bbL^*_{\lambda^+,\kappa}$ let $\Mod^*_\psi$ be
defined as in \ref{b4}(1A),(2),(3) replacing $\uf(\bbB)$ by $\cP(\lambda)$.

\noindent
3) Let $H(\bbL^*_{\lambda^+,\kappa})$ be defined like
 $H(\bbL_{\lambda^+,\kappa}[\bbB])$ in \ref{b4}(5).

\noindent
4) For $\psi \in \bbL^*_{\lambda^+,\kappa}$ let $\spec^*_\psi =
\spec^{1,*}_\psi = \{\|M\|:M \in \Mod^*_\psi\}$; and $\spec^{2,*}_\psi
= \{\|M\|:M \in \Mod^*_\psi$ and $\|M\| = \|M\|^{< \kappa}\}$;
for transparency we will stipulate that from
$\psi$ we can reconstruct $\bbL^*_{\lambda^+,\kappa}$.
\end{definition}

\begin{remark}
\label{b11g}
The following claim essentially tells us that for determining the Hanf
number of $\bbL^{\ba}_{\lambda^+,\kappa}$, we may use the ``worst"
Boolean Algebra, $\bbB^{\fr}_\lambda$ and
$\bbL_{\lambda^+,\kappa}[\bbB^{\fr}_\lambda]$ is essentially equal to
$\bbL^*_{\lambda^+,\kappa}$.  
\end{remark}

\noindent
Parallely to \ref{b5}, \ref{b6}(3):
\begin{claim}  
\label{b10}
1) In the natural definition of $H(\bbL^*_{\lambda^+,\kappa})$ it does
not matter if we use $\spec^{1,*}_\psi$ or $\spec^{2,*}_\psi$ for
$\psi \in \bbL^+_{\lambda^+,\kappa}$.

\noindent
2) For every $\mu < H(\bbL^*_{\lambda^+,\kappa})$ we have $2^\mu <
H(\bbL^*_{\lambda^+,\kappa})$ hence $H(\bbL^*_{\lambda^+,\kappa})$ is
a strong limit cardinal; moreover, of cofinality $> \lambda$.

\noindent
3) $H(\bbL^*_{\lambda^+,\kappa}) < H(\bbL_{2^\lambda)^+,\kappa})$.

\noindent
4) $H(\bbL_{\lambda^+,\kappa}) < H(\bbL^*_{\lambda^+,\kappa}) <
H(\bbL_{(2^\lambda)^+,\kappa})$.

\noindent
5) If $\psi \in \bbL^*_{\lambda^+,\kappa}$ and
$H(\bbL^*_{\lambda^+,\kappa}) \le \sup\{\|M\|:M \in \Mod_\psi\}$,
\then \, $\infty = \sup\{\|M\|:M \in \Mod_\psi\}$.

\noindent
6) For every $\psi_1 \in \bbL^*_{\lambda^+,\kappa}$ there
are $\psi_2,\psi'_2,\psi''_2 \in \bbL^*_{\lambda^+,\kappa}$ such that:
\mn
\begin{enumerate}
\item[$(a)$]   $\spec^*_{\psi_2} = \{\mu:\mu =\mu^{< \kappa} \in
\spec^*_{\psi_1}[\bbB]\} = \spec^{2,*}_{\psi_1}$ 
\sn
\item[$(b)$]   $\spec^*_{\psi'_2} = \{\mu^{< \kappa}:\mu \in 
\spec^{1,*}_{\psi_1}\}$ and
\sn
\item[$(c)$]   $\spec^*_{\psi''_2} = \{\mu:\mu \ge \lambda$ and
$\mu \in \spec^{1,*}_{\psi_1}[\bbB]\}$.
\end{enumerate}
\end{claim}

\begin{PROOF}{\ref{b10}}
Similarly to \ref{b5} and \ref{b6}(3).
\end{PROOF}

\begin{claim}
\label{b12}
1) For every $\psi_1 \in \bbL^*_{\lambda^+,\kappa}$ there is
$\psi_2 \in \bbL_{\lambda^+,\kappa}[\bbB^{\fr}_\lambda]$
   such that $\{\|M\|:M \in \Mod^{\ba}_{\psi_1}\} = \{\|M\|:M \in
\Mod^*_{\psi_2}[\bbB]\}$, that is $\spec^*_{\psi_1} =
\spec_{\psi_2}[\bbB]$.

\noindent
2) For every $\psi_2 \in \bbL_{\lambda^+,\kappa}[\bbB^{\fr}_\lambda]$ there is
$\psi_1 \in \bbL^*_{\lambda^+,\kappa}$ which are as in clause (c).
\end{claim}

\begin{PROOF}{\ref{b12}}
The point is that (A) implies (B) when:
\mn
\begin{enumerate}
\item[(A)]  assume $\bbB$ is the Boolean Algebra generated freely by
  $\langle b_i:i < \lambda\rangle,M$ is a model, $P^M_1 = \{b_i:i <
  \lambda\},P^M_2 = \bbB,Q^M_1 = \cP(\lambda),Q^M_2 = \uf(\bbB),R^M_1
  = \{(c_i,u):u \subseteq \lambda,i \in u\}$ and $R^M_2 = \{((c,D):c
  \in \bbB,D \in \uf(\bbB))$ and $c \in D\},c_{\bar b} \in \bar c(M)$
  and $c^M_b=b$ for $b \in \bbB$
\sn
\item[(B)]  if $N$ is a model of $\Th(M)$ omitting the type $p(x) =
  \{P(x) \wedge x \ne c_b:b \in \bbB\}$ then $(a) \Rightarrow (b)$
  when:
\sn
\begin{enumerate}
\item[(a)]  $N$ satisfies the demands in Definition \ref{b4}(2) of
  $\bbL_{\lambda^+,\kappa}[\bbB^{\fr}_\lambda]$ with $P_2,Q_2,R_2$
  here standing for $P,Q,R$ there
\sn
\item[(b)]  $N$ satisfies the demands in Definition \ref{b9}(1) of
 $L^*_{\lambda^+,\kappa}$ with $P_1,Q_1,R_1$ here standing for $P,Q,R$
 there.
\end{enumerate}
\end{enumerate}
\end{PROOF}

\noindent
Next we have to connect those logics with first order $T$'s.  The easy
part is to start with a Boolean Algebra $\bbB$ and construct a related $T$.
\begin{claim}
\label{b14}
1) For every Boolean Algebra $\bbB$ of cardinality $\le \lambda$ and cardinal
$\kappa \le \lambda^+$ there is $T = T^1_{\bbB,\kappa}$ such that:
\mn
\begin{enumerate}
\item[$(*)_1$]  
\begin{enumerate}
\item[(a)]  $T$ is a first order complete and stable 
\sn
\item[(b)]  $|T| = \lambda$ and $\kappa(T) = \kappa$
\sn
\item[(c)]  $\lambda(T)$ is the cardinality of 
$\uf(\bbB)$, see Definition \ref{a4}(5), \ref{a2}(4), in fact, 
$\bbB_T$ is not much more complicated than $\bbB$
but we shall not elaborate, see \ref{b40} below
\sn
\item[(d)]  $T$ has elimination of quantifiers.
\end{enumerate}
\end{enumerate}
\mn
2) For $\bbB,\lambda,\kappa$ as above there is $T=T^2_{\bbB,\kappa}$
such that:
\mn
\begin{enumerate}
\item[$(*)_2$]  
\begin{enumerate}
\item[${{}}$]  (a),(b) as above
\sn
\item[${{}}$]  $(c) \lambda(T) = \lambda + 2^{\aleph_0}$.
\end{enumerate}
\end{enumerate}
\end{claim}

\begin{PROOF}{\ref{b14}}
Easy, but we elaborate.

\noindent
1) We choose $\tau_*,T_0$ by:
\mn
\begin{enumerate}
\item[$(*)'_1$]
\begin{enumerate}
\item[(a)]  $\tau_* = \tau_{\bbB,\kappa} = \{P_b:b \in
  \bbB\} \cup \{Q_\theta:\theta < \kappa$ is infinite$\} \cup
 \{E_{\theta,i}:\theta < \kappa$ is infinite, $i < \theta\}$ where 
$P_b,Q_\theta$ are unary predicates, $E_{\theta,i}$ 
a binary predicate
\sn
\item[(b)]  universal theory $T_0 \subseteq \bbL(\tau_*)$ is such
  that:  a $\tau_*$-model $M$ satisfied $T_0$ \Iff
\sn
\begin{enumerate}
\item[$(\alpha)$]  $b \mapsto P^M_b$ embeds $\bbB$ into the
  Boolean Algebra  $\cP(P^M_{1_{\bbB}})$ so $P^M_{0_{\bbB}} = \emptyset$
\sn
\item[$(\beta)$]  $\langle P^M_{1_{\bbB}} \rangle \char 94
\langle Q^M_\theta:\theta < \kappa\rangle$ are pairwise disjoint
\sn
\item[$(\gamma)$]  $E^M_{\theta,i}$ is an equivalence
  relation on $Q^M_\theta$ so $a E^M_{\theta,i} b \Rightarrow a,b \in
 Q^M_\theta$
\sn
\item[$(\varepsilon)$]   if $i<j < \theta$ then
  $E^M_{\theta,j}$ refines $E^M_{\theta,i}$.
\end{enumerate}
\end{enumerate}
\end{enumerate}
\mn
So
\mn
\begin{enumerate}
\item[$\oplus_1$]
\begin{enumerate}
\item[(a)]  $T_0 \subseteq \bbL(\tau_*)$ is a well defined universal theory
\sn
\item[(b)]  $\Mod_{T_0}$ has amalgamation and the $\JEP$.
\end{enumerate}
\end{enumerate}
\mn
Let
\mn
\begin{enumerate}
\item[$\oplus_2$]  $\bbT$ is the set of $\tau \subseteq \tau_*$
  satisfying:
\sn
\begin{enumerate}
\item[(a)]  $P,P_{1_{\bbB}},P_{0_B} \in \tau$
\sn
\item[(b)]  $E_{\theta,i} \in \tau \Rightarrow Q_\theta \in \tau$
\sn
\item[(c)]  if $\bbB \models ``b \cap c = a \wedge -b = d"$ then
  $\{P_b,P_c\} \subseteq \tau_1 \Rightarrow \{P_a,P_d\} \subseteq
  \tau$
\end{enumerate}
\sn
\item[$\oplus_3$]  for $\tau \in \bbT$ let $T_{0,\tau}$ be defined like
  $T_0$ but restricting ourselves to predicates from $\tau$.
\end{enumerate}
\mn
Now
\mn
\begin{enumerate}
\item[$\oplus_4$]  for $\tau \in \bbT$
\sn
\begin{enumerate}
\item[(a)]  if $M$ is a $\tau$-model of $T_{0,\tau}$, \then \, $M$ can
  be expanded to a $\tau_*$-model of $T_0$
\sn
\item[(b)]  $T_{0,\tau}$ has the $\JEP$
\sn
\item[(c)]  $T_{0,\tau}$ has the amalgamation property
\sn
\item[(d)]  if $M_1 \subseteq M_2$ are models or $T_{0,\tau}$ and
  $\tau \subseteq \tau_1 \in \bbT$ and $N_1$ is a $\tau_1$-model
  expanding $M_2$ \then \, there is a $\tau_1$-model $N_2$ expanding
$M_1$ and extending $N_1$.
\end{enumerate}
\end{enumerate}
\mn
[Why?  Easy, e.g. clause (b) by disjoint union.]
\mn
\begin{enumerate}
\item[$\oplus_5$]  For finite $\tau \in \bbT,T_{0,\tau}$ has a model
  completion called $T_{1,\tau}$ which has elimination of quantifiers.
\end{enumerate}
\mn
[Why?  Because $\tau$ is a relational finite vocabulary and
$T_{0,\tau}$ is univesal with $\JEP$ and amalgamation.]
\mn
\begin{enumerate}
\item[$\oplus_6$]   If $\tau_1 \subseteq \tau_2$ are from $\bbT$ \then
  \, $T_{1,\tau_1} \subseteq T_{1,\tau}$.
\end{enumerate}
\mn
[Why?  By $\oplus_4(d) + \oplus_5$.]
\mn
\begin{enumerate}
\item[$\oplus_7$]   $T = T^1_{\bbB,\kappa} := 
\cup\{T_{1,\tau}:\tau \in \bbT$ finite$\}$
is the model completion of $T_0$ and has elimination of quantifiers.
\end{enumerate}
\mn
[Why?  Follows from the above.]
\mn
\begin{enumerate}
\item[$\oplus_8$]
\sn
\begin{enumerate}
\item[(a)]  If $\tau \in \bbT$ is finite, \then \, $T_{1,\tau}$ is
  $\aleph_0$-categorical and $\aleph_0$-stable
\sn
\item[(b)]  $T$ is stable
\sn
\item[(c)]  $\kappa(T) = \kappa$
\sn
\item[(d)]  $|\lambda'(T)| = |\bbB| + \aleph_0$
\sn
\item[(e)]  $\lambda(T) = \min\{\mu:\mu \ge \lambda$ and 
$\mu^{< \kappa} = \mu\}$.
\end{enumerate}
\end{enumerate}
\mn
[Why?  Consider the monster $\gC = \gC_{T_{1,\tau}}$ and use
automorphisms.]

So $T = T^1_{\bbB,\kappa}$ from $\oplus_7$ is as promised.

\noindent
2) We use $T_0$ such that $(*)'_2$ below holds and continue as above.
\mn
\begin{enumerate}
\item[$(*)'_2$]  as in $(*)'_1$ above but
\sn
\begin{enumerate}
\item[$(a)$]   we add $Q_0,E_{0,n}(n < \omega)$ wih $Q_0$ unary and
  $E_{0,n}$ binary
\sn
\item[$(b)$] 
\sn
\begin{enumerate}
\item[$(\beta)$]  also $Q^M_0$ is disjoint to $Q^M_\theta(\theta \in
  [\aleph_0,\kappa))$ and to $P^M_{1_{\bbB}}$
\sn
\item[$(\zeta)$]  $E^M_{0,n}$ is an equivalence relation on $P^M_0$
\sn
\item[$(\eta)$]  $E^M_{0,0}$ has one equivalence class
\sn
\item[$(\theta)$]  $E^M_{0,n+1}$ refines $E^M_{0,n}$ and
  divides each $E^M_{0,n}$ equivalence class to at most 2.
\end{enumerate}
\end{enumerate}
\end{enumerate}
\end{PROOF}

\begin{discussion}
\label{b16}
1) We like to translate ``$M \models ``\psi,\psi \in
\bbL_{\lambda^+,\kappa}"$ to ``$M \in \Mod_{\bold t}$", that is, when
$\kappa(T) \ge \kappa$ and, in particular, when $\kappa > \aleph_0$.  However, 
the following is the ``translation of $\psi \in
\bbL_{\lambda^+,\kappa}(\tau_0)$"; i.e. it deals strictly with the logic
$\bbL_{\lambda^+,\kappa}$; in particular a Boolean Algebra $\bbB$ is not
present.  Our aim is to do some of the work of \ref{b24} in which we are really
interested.  So \ref{b20} is not directly related to $\bold t$'s! as there
is no saturation requirement; moreover stability appears neither in
\ref{b20} nor in \ref{b24}. 

\noindent
2) Note that in \ref{b20} we can let $\kappa_1$ be such that $\kappa =
\kappa^+_1$ or $\kappa_1 = \kappa$ is a limit cardinal and let
$\Upsilon = \kappa_1 +1$ and omit $F_{\kappa_1},P_{\kappa_1}$.
\end{discussion}

\begin{theorem}  
\label{b20}
\underline{The $\bbL_{\lambda^+,\kappa}$-representation Theorem}

Assume $\psi \in \bbL_{\lambda^+,\kappa}(\tau_0)$, so of course, $|\tau_0| \le
\lambda$.  Let $\Upsilon$ be $\kappa$ if $\kappa \le \lambda$ and
$\lambda +1$ if $\kappa = \lambda^+$.

 \Then \, we can find a tuple $(\tau_1,T_1,p(x),\bar F,\bar P)$ 
such that (for $\bar F,\bar P$ as below):
\mn
\begin{enumerate}
\item[(A)]
\begin{enumerate}
\item[(a)]  $\tau_1$ is a vocabulary $\supseteq \tau_0$ of
cardinality $\lambda$
\sn
\item[(b)]  $\bar F$ is a sequence of unary function symbols with no
repetitions of length $\Upsilon$, new (i.e. from $\tau_1 \backslash
\tau_0$), let $\bar F = \langle F_i:i < \Upsilon\rangle$
\sn
\item[(c)]  $\bar P$ is a sequence of unary predicates with
no repetitions of length $\Upsilon$, new 
(i.e. from $\tau_1 \backslash \tau_0$), 
let $\bar P = \langle P_i:i < \Upsilon\rangle$
\sn
\item[(d)]  $T_1$ is a first order theory in the vocabulary $\tau_1$
\sn
\item[(e)]  $p(x)$ is $\{P_*(x) \wedge x \ne c_i:i < \lambda\}$,
an $\bbL(\tau_1)$-type (even quantifier-free), 
so $P_*$ is a unary predicate and $c_i$ for $i < \lambda$ individual
constants, all new
\end{enumerate}
\sn
\item[(B)]  the following conditions on a $\tau_0$-model $M_0$ are
equivalent
\sn
\begin{enumerate}
\item[(a)]  $M_0 \models \psi$ and $\|M_0\| = \|M_0\|^{< \kappa} +
  \lambda^{< \kappa}$
\sn
\item[(b)]  there is a $\tau_1$-expansion $M_1$ of
$M_0$ to a model of $T_1$ omitting $p(x)$ such that:
\sn
\begin{enumerate}
\item[$(\alpha)$]  $\langle P^{M_1}_i:i < 
\Upsilon\rangle$ is a partition of $|M_1|$
\sn
\item[$(\beta)$]  if $i < \Upsilon$ and $a_j \in M_1$ for
$j<i$ then for some $b \in P^{M_1}_i$ we have 
$j<i \Rightarrow F^{M_1}_j(b) = a_j$.
\end{enumerate}
\end{enumerate}
\end{enumerate}
\end{theorem}

\begin{PROOF}{\ref{b20}}  
Note that as $\psi$ has no free variables, \wilog \, every subformula
$\varphi$ of $\psi$ has a set of free variables equal to $\{x_i:i <
\varepsilon\}$ for some $\varepsilon = \varp_\varphi < \kappa$ such
that if $\varphi$ is a subformula of $\psi$ and $\varphi =
\bigwedge\limits_{i < j} \varphi_i$ then $\varp_{\varphi_i} =
\varp_\varphi$. 

Let $\Delta$ be the set of subformulas of $\psi$ so \wilog \, (a
syntactial rewriting) there is a list
$\langle \varphi_i(\bar x_{[\varepsilon(i)]}):i < i(*)\rangle$ for
some $i(*) \le \lambda$ of $\Delta$ such that $\varepsilon(0) = 0,\varphi_0 =
\psi$ and $\bar x_{[\varepsilon(i)]}$ is a sequence of length $< \kappa$
of variables, in fact, $\bar x_{[\varepsilon(i)]} = \langle
x_\varepsilon:\varepsilon < \varepsilon(i)\rangle$ and $\varepsilon(i) <
\kappa$.

For any $\tau_0$-model $M$ such that $\|M\| = \|M\|^{< \kappa} +
\lambda^{< \kappa}$, we say $N$ codes $M$ when:
\mn
\begin{enumerate}
\item[$(*)$]
\begin{enumerate}
\item[(a)]  $N$ expands $M$
\sn
\item[(b)]  $\langle F^N_i:i < \Upsilon\rangle,\langle P^N_i:i <
\Upsilon \rangle$ satisfies $(B)(b)(\alpha),(\beta)$ of the theorem
(with $N$ instead of $M_1$)
\sn
\item[(c)]  $Q^N_i = \{b \in P^N_{\varepsilon(i)}:M \models
\varphi_i[\langle F_\varepsilon(b):\varepsilon <
  \varepsilon(i)\rangle]\}$ for $i<i(*)$
\sn
\item[(d)]  $\langle c^N_i:i < \lambda\rangle$ are pairwise distinct
and $P^N_* = \{c^N_i:i < \lambda\}$
\sn
\item[(e)]  if $\varphi_i(\bar x_{\varepsilon(i)}) =
\bigwedge\limits_{j < j(i)} \varphi_{i,j}(\bar x_{\varepsilon(i)})$ 
so for some $\bfi(i,j) < i(*)$ we have
$\varphi_{i,j}(\bar x_{\varepsilon(i)}) = \varphi_{\bold i(i,j)}(\bar
x_{\varepsilon(\bold i(i,j))})$ and
so $\varepsilon(\bold i(i,j)) =
\varepsilon(i)$ \then \, $F_{1,i} \in \tau(N)$ is unary and for 
$b \in P^N_{\varepsilon(i)}$ we have:
\sn
\begin{enumerate}
\item[$(\alpha)$]  $N \models ``F_{1,i}(b) = c_j \wedge
\neg \varphi_i(\langle F_\varepsilon(b):\varepsilon < 
\varepsilon(i)\rangle)"$ implies $M \models 
\neg \varphi_{i,j}(\langle F_\varepsilon(b):\varepsilon <
\varepsilon(i)\rangle)$ which means: 
if $\varphi_{i,j} = \varphi_{\bold i(i,j)}$ and $N \models
``\neg Q_i(b) \wedge c_j = F_{1,i}(b)"$ then 
$M \models ``\neg Q_{\bold i(i,j)}[b]"$ 
and, of course
\sn
\item[$(\beta)$]   if $M \models
  \varphi_i(\langle f_\varepsilon(b):\varepsilon <
  \varepsilon(i)\rangle)$ and $j < \varepsilon(i)$ then $M \models
\varphi_{i,j}(\langle F_\varepsilon(b):\varepsilon < \varepsilon(i)\rangle)$
\end{enumerate}
\sn
\item[(f)]   if $\varphi_i(\bar x_{\varepsilon(i)}) =
(\exists \bar x_{[\varepsilon(i),\zeta(i))}) \varphi_{j_1(i)}(\bar
  x_{\varepsilon(i)},\bar x_{[\varepsilon(i),\zeta(i))})$ and
  $F_\varepsilon(b) = a_\varepsilon$ for
$\varepsilon < \varepsilon(i)$ \then \, $(\alpha)
\Leftrightarrow (\beta)$ where
\sn
\begin{enumerate}
\item[$(\alpha)$]  $M_1 \models 
\varphi_i[\langle a_\varepsilon:\varepsilon < \varepsilon(i)\rangle]$
equivalently $M_1 \models \varphi_1[\langle F_\varp(b):\varp <
\varp(i)\rangle]$ 
\sn
\item[$(\beta)$]  $M_1 \models (\exists y)
\varphi_{j_1(i)}(\langle a_\varepsilon:
\varepsilon < \varepsilon (i)\rangle,\langle F_\zeta(y):
\zeta \in [\varepsilon(i),\zeta(i)]\rangle$.
\end{enumerate}
\end{enumerate}
\end{enumerate}
\mn
Now let
\mn
\begin{enumerate}
\item[$\boxplus$]  
\begin{enumerate}
\item[(a)]  $\tau_1$ is $\tau_\psi \cup
  \{F_\varepsilon,P_\varepsilon:\varepsilon < \Upsilon\} \cup \{Q_i:i <
  i(*)\} \cup \{F_{1,i}:i <i(*)$ and $\varphi_i$ 
 is a conjunction$\}$
\sn
\item[(b)]  $T_1 = \cap\{\Th(N)$: there is $M$, a $\tau_0$-model of
$\psi$ such that $\|M\| = \|M\|^{< \kappa} + \lambda$ and 
$N$ code $M\}$ 
\sn
\item[(c)]  $p(x) = \{P_*(x) \wedge x \ne c_i:i < \lambda\}$.
\end{enumerate}
\end{enumerate}
\mn
Now check that
\mn
\begin{enumerate}
\item[$\oplus$]  $(\tau_1,T_1,p(x),\bar F,\bar P)$ is as required.
\end{enumerate}
\end{PROOF}

\begin{remark}
\label{b21g}
So how does \ref{b20} help for our main aim?  It starts to translate
$\psi \in \bbL_{\lambda^+,\kappa}(\tau_0)$ to $\bold t = (\tau_1,T_1,p(x))$, so
instead having blocks of quantifiers $(\exists \bar
x_{[\varepsilon]}),\varepsilon < \kappa$ we have $(\exists x)$, i.e. 
by the sequence of functions $\langle F_i:i < \varepsilon\rangle$ we code any 
$\varepsilon$-tuple by one element.

This will help later to make ``the $\tau(T_{\bold t})$-reduct is saturated"
equivalent to the existence of suitable coding.
\end{remark}

\noindent
Recalling Definition \ref{b4}(6) of $\bbL_{\lambda^+,\kappa}[\bbB]$,
we get the section main result: translating from $\psi \in
\bbL_{\lambda^+,\kappa}[\bbB]$ to a representation, naturally more
complicated than the one for $\psi \in \bbL_{\lambda^+,\aleph_0}$.
\begin{theorem}  
\label{b24}
\underline{The $\bbL_{\lambda^+,\kappa}[\bbB]$-representation theory}

Assume $\bbB$ is a Boolean Algebra of cardinality $\le \lambda$ and
for notational transparency $b \in \bbB \cap \alpha < \lambda
\Rightarrow b \ne \alpha$ and $\psi \in
\bbL_{\lambda^+,\kappa}[\bbB](\tau_0)$.  \Then \, we can find a tuple 
$(\tau_1,T_1,p(x),\bar F,\bar P)$ 
such that (for $\bar F,\bar P$ as below):
\mn
\begin{enumerate}
\item[(A)]  as in \ref{b20}
\sn
\item[(B)]  the following conditions on a $\tau_0$-model $M_0$ are
equivalent:
\sn
\begin{enumerate}
\item[(a)]  $M_0 \in \Mod^2_\psi[\bbB]$, so
$M_0 \models \psi$ and $\|M_0\| = \|M_0\|^{< \kappa} + \lambda^{< \kappa}$
\sn
\item[(b)]  there is a $\tau_1$-expansion $M_1$ of
$M_0$ to a model of $T_1$ omitting $p(x)$ such that:
\sn
\begin{enumerate}
\item[$(\alpha)$]  $\langle P^{M_1}_i:i < 
\Upsilon\rangle$ is a partition of $|M_1|$
\sn
\item[$(\beta)$]  if $i < \Upsilon$ and $a_j \in M_1$ for
$j<i$ then for some $b \in P^{M_1}_i$ we have 
$j<i \Rightarrow F^{M_1}_j(b) = a_j$
\sn
\item[$(\gamma)$]  $c_b(b \in \bbB)$ are 
individual constants (in $\tau_1 \backslash \tau_0$) with no repetition, 
 $P,Q \in \tau_1$ unary, $R \in \tau_1$ binary 
\sn
\item[$(\delta)$]  $P^{M_1} = \{c^{M_1}_b:b \in \bbB\}$
\sn
\item[$(\varepsilon)$]  $R^{M_1} \subseteq P^{M_1} \times Q^{M_1}$
\sn
\item[$(\zeta)$]   for every $b \in Q^{M_1}$ the set
$u(b,M_1) := \{c_b \in P^{M_1}:(c_b,b) \in R^{M_1}\}$ 
 is an ultrafilter of $\bbB$
\sn
\item[$(\eta)$]  for every ultrafilter $D$ of the
Boolean Algebra $\bbB$ there is one
and only one $b \in Q^{M_1}$ such that $u(b,M_1) = D$.
\end{enumerate}
\end{enumerate}
\end{enumerate}
\end{theorem}

\begin{PROOF}{\ref{b24}}  
First, note that $P,Q,c_b(b \in \bbB)$ are in $\tau_\psi$ as in
Definition \ref{b4}.  Second, we repeat the proof of \ref{b20} or
just quote it:
\mn
\begin{enumerate}
\item[$(*)_1$]  there is $\tau_* \supset \tau_\psi,|\tau_*| = \lambda$
  with $F_\varp,P_\varp,F_{1,\varp},c_\varp,Q \in \tau_*$ as there,
  i.e. satisfying clauses (A)(a)-(e).
\end{enumerate}
\mn
Third, we prove clause (B) of \ref{b24}.  The direction (B)(b)
$\Rightarrow$ (B)(a) holds as in \ref{b20}.  For the other direction,
assume $M_0 \in \Mod^0_\psi[\bbB]$ and we choose $M_1$ as in
\ref{b20}$(B)(\alpha),(\beta)$.

Lastly, clauses (B)(b)$(\gamma)-(\eta)$ holds because $\psi \in
\bbL_{\lambda^+,\kappa}[\bbB]$ and $M_1$ expands $M_0$.
\end{PROOF}

\begin{remark}  
\label{b26}
1) The only non-``$\bbL_{\lambda^+,\kappa}$ demand" in clause (B) of
\ref{b24} is in $(b)(\eta)$, the existence, 
this is not expressible by a sentence of
$\bbL_{\lambda^+,\kappa}$, even with extra predicates.

\noindent
2) As indicated above, $\bbB^{\fr}_\lambda$ is the ``worst, most
complicated Boolean Algebra" for our purpose.  So it is natural to
wonder about the order among the relevant Boolean Algebras, so
\ref{b40}, \ref{b43} try to deal with it.
\end{remark}

\begin{definition}
\label{b40}
1) We define a two-place relation $\le^*_{\lambda^+,\kappa}$ among the
Boolean Algebras $\bbB$ of cardinality $\le \lambda$

$\bbB_1 \le^*_{\lambda^+,\kappa} \bbB_2$ \Iff \,: there is a sentence
$\psi_2 \in \bbL_{\lambda^+,\kappa}[\bbB_2]$, unary predicates $P_1,Q_1
\in \tau_\psi$ and binary predicate $R_2$ and individual constants
$c^1_b(b \in \bbB_1)$ from $\tau_\psi$ such that:
\mn
\begin{itemize}
\item  if $M \models \psi_2$ then $P^M_1 = \{(c^1_b)^M:b \in
  \bbB_1\}$ and $R^M_1 \subseteq P^M_1 \times Q^M_1$ and $\langle
  (c^1_b)^M:b \in \bbB_1\rangle$  satisfies the demands
  in \ref{b4}(2).
\end{itemize}
\mn
2) We let $\equiv^*_{\lambda^+,\kappa}$ be defined by $\bbB_1
\equiv^*_{\lambda^+},\kappa \bbB_2$ \Iff \, $\bbB_1
\le^*_{\lambda^+,\kappa} \bbB_2$ and $\bbB_2 \le^*_{\lambda^+,\kappa} \bbB_1$.
\end{definition}

\begin{claim}
\label{b43}
1) $\le^*_{\lambda^+,\kappa}$ is a quasi-order on the class of Boolean
Algebras of cardinality $\le \lambda$.

\noindent
2) Hence $\equiv^*_{\lambda,\kappa}$ is an equivalence relation with
being isomorphic refining it.

\noindent
3) In \ref{b14}(1) we have $\bbB_T \equiv^*_{\lambda^+,\kappa} \bbB$
where $T=T^1_{\bbB,\lambda}$.

\noindent
4) If $\bbB_1 \le^*_{\lambda,\kappa} \bbB_2$ \then \, for every
$\psi_1 \in \bbL_{\lambda^+,\kappa}[\bbB_1]$ there is $\psi_2 \in
\bbL_{\lambda^+,\kappa}[\bbB_2]$ such that:
\mn
\begin{enumerate}
\item[$(a)$]  $\spec^2_{\psi_1} = \spec^2_{\psi_2}$ 
\sn
\item[$(b)$]  if $M_1$ is a $\tau(\psi_1)$-model then $M_1 \in
  \Mod^2_{\psi_1}$ \Iff \, $M_1 = M_2 \rest \tau_{\psi_1}$ for some
  $M_2 \in \Mod^2_{\psi_2}$; pedantically we should have an embedding
  $\pi$ of $\tau_{\psi_1}$ into $\tau_{\psi_2}$ and demand $M_1 = (M_2
  \rest \Rang(\pi))^{[\pi]}$, naturally defined.
\end{enumerate}
\mn
5) If $\bbB$ is a Boolean Algebra of cardinality $\le \lambda$ then
$\bbB \le^*_{\lambda,\kappa} \bbB^{\fr}_\lambda$.
\end{claim}

\begin{PROOF}{\ref{b43}}
1) Easy but we elaborate; so asume $\bbB_1,\bbB_2,\bbB_3$ are Boolean
Algebras of cardinality $\le \lambda$ and $y_1
\le^*_{\lambda^+,\kappa} \bbB_2$ and $\bbB_2 \le^*_{\lambda^+,\kappa}
\bbB_3$.  Hence for $\ell=1,2$ ther is a sentence $\psi_\ell$ and
$P_1,Q-1,R_2 c^1_b (b \in \bbB_\ell)$ from $\tau_{\psi_\ell}$
witnessing it, and let $P,Q,R,c_b(b \in \bbB_{\ell +1})$ be as
promised in Definition \ref{b4} for $\psi \in
\bbL_{\lambda^+,\kappa}[\bbB_{\ell +1}]$.  We can find disjoint
vocabularies $\tau_1,\tau_2$ and function $h_1,h_2$ such that:
\mn
\begin{enumerate}
\item[$(*)$]  for $\ell=1,2$ the function $h_\ell$ is a one-to-one
  functino from $\tau(\psi_\ell)$ onto $\tau_\ell$, preserving ``being
  a predicate/function symbol/individual constant" and preserving the
  arity; let $\varphi_\ell$ be the image of $\psi_\ell$ under $h$.
\end{enumerate}
\mn
Lastly, let $\varphi$ be the conjunction of:
\mn
\begin{enumerate}
\item[(a)]  $\varphi_1,\varphi_2$
\sn
\item[(b)]  $h_1(P),h_2(P_1)$ are equivalent $(\forall x)[(h_1(P)(x)
  \equiv (h_2(P_1))(x)]$
\sn
\item[(c)]  also $h_1(Q),h_1(R),h_1(c_b)$ are equivalent to
  $h_2(Q_2),h_2(R_2),h_2(C^1_b)$ respectively.
\end{enumerate}
\mn
The rest should be clear.

\noindent
2) Follows from part (1).

\noindent
3) Let us fix $m \ge 1$ and we shall analyze $\bbB_{T_1,m}$.  Let
$\Lambda_1 = \{\eta:\eta$ is a sequence of length $m$ with range
included in $\Theta\}$ where $\Theta = \{\theta:\theta < \kappa$
infinite$\} \cup \{0\}$.

For $\theta \in \Theta$ let $\varphi_0(x) = Q_\theta(x)$, interpreting
$Q_0$ as $P$.  Next let $\Lambda_0 = \{\eta \rest u:u \subseteq m$ and
$\eta \in \Lambda_1\}$ and for $\eta \in \Lambda_0$ let
$\varphi_{\eta,\ell}(\bar x_{[m]}) = \bigwedge\limits_{\ell < m}
\varphi_{\eta(\ell)}(x_\ell)$ and for $\nu \in \Lambda_0 \backslash
\Lambda_1$ let $\Lambda_{0,\nu} = \{\eta \in \Lambda_0:\nu \subseteq
\eta$ and $\Rang(\eta) \backslash \Rang(\nu)$ is a singleton?$\}$.

Lastly
\mn
\begin{enumerate}
\item[$(*)$]  $\bbB_{T,m,\eta} = \bbB_{T,m} \rest \{\bar a:a \le
  \varphi_\eta(\bar x_{[m]})/\equiv_T\}$ for $\eta \in \Lambda_0$
\sn
\item[$(*)$]  if $\nu \in \Lambda_0 \backslash \Lambda_1$, \then \, 
$\bbB_{T,m,\nu}$ is the direct sum of $\langle \bbB_{T,m,\eta}:\eta
\in \Lambda_{0,\nu}\rangle$
\sn
\item[$(*)$]  $\bbB_{T,m,\emptyset} = \bbB_{T,m}$
\sn
\item[$(*)$]  if $D \in \uf(\bbB_{T,m})$ \then \,
\sn
\begin{enumerate}
\item[$\bullet_1$]  for some $\eta \in \Lambda_0,\varphi_\eta(\bar
  x_{[m]})/\equiv_T \in D,|\dom(\eta)|$ minimal
\sn
\item[$\bullet_2$]  so $D$ is deteremined by $\eta$ and $D \rest
  \bbB_{T,m,\eta} \in \uf(\bbB_{T,m,n})$
\sn
\item[$\bullet_3$]  if $\eta \in \Lambda$ then $\eta$ determines $D$
\sn
\item[$\bullet_4$]  if $\Rang(\eta)$ is minimal it is a singleton, so
\end{enumerate}
\sn
\item[$(*)$]  above if $\Rang(\eta) = \{\theta\},\theta \ge \aleph_0$
  \then \, $\bbB_{T,m,\eta}$ is isomorphic to $\bbB_{T^{\eq}_\theta}$,
  i.e. $\theta$-equivalence relation (see Definition \ref{a16})
\sn
\item[$(*)$]  above if $\Rang(\eta) = \{0\}$, then $\bbB_{T,m,\eta}$
  is isomorphic to the direct sum of $|\{(e,a):e$ an equivalence
  relation on $\dom(\eta)\}|$ and $a$ is an $e$-equivalence class
\sn
\item[$(*)$]  
\begin{enumerate}
\item[(a)]  the number of ultrafilters on $\bbB_{T,m} \, (m \ge 0)$ is
  $|\uf(\bbB)|$ if $|\uf(\bbB)| \ge \lambda$
\sn
\item[(b)]   $|\uf(\bbB_{T,M})| = \sup\{\theta:\theta < \kappa\} +
  \aleph_0$.
\end{enumerate}
\sn
\item[$(*)$]  $\lambda(T) = |\uf(\bbB)| + \sup\{\theta:\theta <
  \kappa\} + \aleph_0$.
\end{enumerate}
\mn
Also the $\equiv^*_T$ is easy.

\noindent
4) Read the definition.

\noindent
5) Holds by \ref{b6}(1)(a).
\end{PROOF}
\newpage

\section {Real equality for each $T$}
\bigskip

\subsection {Answering the Original Question and the New One}\
\bigskip

The original question for this work 
was about the strictly stable case, i.e. fixing
$\kappa > \aleph_0$, dealing with $\{\bold t \in \bold
N_\lambda:\kappa(T_{\bold t}) = \kappa\}$, so we deal with this case
first.

In this case Theorem \ref{c2} tells us that for strictly stable $T$
and $\lambda \ge |T|$, the family of classes $\Mod_{\bold t}$ for
$\bold t \in \bold N_{\lambda,T}$ and the family of classes
$\Mod^2_\psi[\bbB]$ for $\psi \in \bbL_{\lambda^+,\kappa}[\bbB]$ where
$\kappa = \kappa_r(T)$ and $\bbB$ is the Boolean algebra $\bbB_T$ from
\ref{a4}(2),(2A),(3),(3A) are very similar.  
How this is proved?  For one direction, we start with
$\bold t \in N_{\lambda,T}$; so the (essential)
non-first order part of the demand
$M \in \Mod_{\bold t}$ is ``$M \rest \tau(T_{\bold t})$ is saturated".
At first glance we need (in addition to the first order theory and the omission
of a type) to say some things on eliminating $u \in [M]^{<\|M\|}$ and
relation on it, but because of $T$ being stable it can be (see
\ref{a9}) expressed by the equivalence of:
\mn
\begin{enumerate}
\item[$(a)$]  $M \rest \tau(T_{\bold t})$ is $\kappa_r(T)$-saturated
\sn
\item[$(b)$]  if $\bold I \subseteq M$ is an infinite indiscernible
  set in $M \rest \tau(T_{\bold t}),|\bold I| = \aleph_0$ then we can
  find an indiscernible set $\bold J \supseteq \bold I$ in $M \rest
  \tau(T_{\bold t})$ of cardinality $\|M\|$.
\end{enumerate}
\mn
So the use of $\bbL_{\lambda^+,\kappa}$ where $\kappa = \kappa_r(T)$
is natural.  If $2^{|T|} \le \lambda$ this is obvious but otherwise we
have to be more careful.  We use the Boolean algebra $\bbB = \bbB_T$
and the use of $\psi \in \bbL_{\lambda^+,\kappa}[\bbB]$ rather than
$\bbL_{\lambda^+,\kappa}$ to express $M \rest \tau(T_{\bold t})$ is
$\aleph_0$-saturated, so by $\kappa_r(T)$-sequence
homogeneity this is enough.

Note that on the one hand $M \in \Mod_{\bold t} \Rightarrow \|M\| \in
\bold C_T = \{\mu:\mu = \mu^{< \kappa(T)} + \lambda(T)\}$, see
\ref{b3} but on the
other hand for $\psi \in \bbL_{\lambda^+,\kappa}[\bbB],M \models \psi$
does not imply it.  Still we know that $\spec^1_\psi = \{\|M\|:M \models
\psi\}$ and $\spec^2_\psi = \spec^1_\psi \cap \bold C_T$ are closed
enough, see Claim \ref{b5}, in particular \ref{b5}(1).  Recall that $\bbB =
\bbB^{\fr}_\lambda$ is the worst case.

For superstable $T$ (for the case we fix $(\lambda,T))$, the case, of
e.g. $ = \Th({}^\omega 2,E_n)_n,E_n = \{(\eta,\nu):\eta,\nu \in
{}^\omega 2,\eta \rest n = \nu \rest n\}$ makes us work somewhat more.

\begin{theorem}  
\label{c2}
Assume $T$ is a stable first order complete of cardinality $\le
\lambda$ and $\kappa = \kappa_r(T)
= \min\{\theta:\theta$ regular and $\theta \ge \kappa(T)\}$ and
$\lambda(T) = \min\{\lambda:T$ stable in $\lambda\}$, see
\ref{b13}(3), and let $\bbB = \bbB_T$, see Definition \ref{a4}(3A).

Assume further that $\kappa(T) > \aleph_0$ (i.e. $T$ is not
superstable).

\noindent
1) We have $\{\spec_{\bold t}:\bold t \in \bold N_{\lambda,T}\} 
= \{\spec^2_\psi[\bbB]:\psi \in \bbL_{\lambda^+,\kappa}[\bbB]\}$.  

\noindent
2) If $\tau_0 = \tau_T$ and $\psi_0 = \wedge\{\varphi:\varphi \in
 T\}$ \underline{or} just $\tau_T \subseteq \tau_0,|\tau_0| \le
   \lambda,\psi_0 \in \bbL_{\lambda^+,\kappa}[\bbB](\tau_0)$ and $M
   \in \Mod_{\psi_0}[\bbB] \Rightarrow M \models T$ \then \, there is
$\bold t \in \bold N_{\lambda,T}$ such that $\spec^2_{\psi_0}[\bbB] 
= \spec_{\bold t}$.

\noindent
3) If $\bold t \in \bold N_{\lambda,T}$ \then \, for some $\psi_1 \in
\bbL_{\lambda^+,\kappa}[\bbB](\tau_1),\tau_1 \supseteq
   \tau(T_2)$ and $\spec^1_{\psi_1}[\bbB] = \spec_{\bold t} =
   \spec^2_{\psi_1}[\bbB]$.
\end{theorem}

\begin{remark}  
\label{c4}
The proof gives more: that the two contexts have the same PC classes.
This proof is divided to two subsections each to one direction.
\end{remark}

\begin{PROOF}{\ref{c2}}
1) By parts (2),(3).

\noindent
2) By \S(2C) below.

\noindent
3) By \S(2B) below, i.e. by \ref{d4} noting \ref{d2}.
\end{PROOF}

\begin{conclusion}
\label{c6}
If $T$ is first order complete stable theory, $\kappa = \kappa(T)$ and
$|T| \le \lambda$ \then \, $H(\bold N_{\lambda,T})$ is bigger than
$H(\bbL_{\lambda^+,\kappa})$ but smaller than $H(\bbL_{(2^\lambda)^+,\kappa})$.
\end{conclusion}

\begin{PROOF}{\ref{c6}}
First assume $T$ is strictly stable, i.e. $\kappa(T) > \aleph_0$.
The ``bigger than $H(\bbL_{\lambda^+,\kappa})$" follows by \ref{c2}(2)
recalling \ref{b5}(4), the first inequality. The ``smaller than
$H(\bbL_{(2^\lambda)^+,\kappa})$" follows by \ref{c2}(3) recalling
\ref{b5}(4), the second and third inequality.   We are left with the
case $T$ is superstable, but then we quote \cite[Th.1.2]{BlSh:992}, or
see \ref{d8}, \ref{d12} below.
\end{PROOF}
\bigskip

\subsection {Given $\bold t \in \bold N_{\lambda,1}$} \
\bigskip

\begin{hypothesis}  
\label{d2}
For this subsection we are given $\bold t = (T,T_1,p) \in
\bold N_{\lambda,T}$ such that $T$ is complete first order stable
 so $\lambda \ge |T_1| \ge |T|$ and let $\bbB = \bbB_T,\kappa = 
\kappa_\tau(T)$;  \wilog \,:
\mn
\begin{enumerate}
\item[$(a)$]   $P,Q,R,c_b(b \in \bbB)$ are \underline{not} in
  $\tau(T_1)$ and with no repetition
\sn
\item[$(b)$]   $P,Q$ are unary predicates, $R$ is a binary predicate,
  $c_b$ individual constants
\sn
\item[$(c)$]   $\tau_2 = \tau(T_1) \cup \{P,Q,R,c_b:b \in \bbB\}$.
\end{enumerate}
\end{hypothesis}

\begin{claim}  
\label{d4}
Assume $\kappa > \aleph_0$.  
There is $\psi \in \bbL_{\lambda^+,\kappa}[\bbB](\tau_1)$ 
such that $\Mod_{\bold t} = \{N \rest \tau(T_1):
N \models \psi$ so $\tau(N) = \tau(\psi)\supseteq \tau_1\}$.
\end{claim}

\begin{PROOF}{\ref{d4}}  
Note that below proving \ref{d8}, \ref{d12} we use this proof stating
the changes; there $\kappa(T) = \aleph_0$, i.e. $T$ is superstable.
\medskip

\noindent
\underline{Stage A}:

\Wilog \, we can replace $T$ by $T^{\eq}$ (no need for new elements:
we can extend $T_1$ to have a copy of $M^{\eq}$ with new predicates and an
isomorphism).  The use of $T^{\eq}$ is anyhow just for transparency.
For $\theta =
\cf(\theta) < \kappa_r(T)$ choose a sequence $\bar\varphi_\theta =
\langle \varphi_{\theta,i}(x,\bar y_{\theta,i}):i < \theta\rangle$
witnessing $\theta < \kappa_r(T)$ equivalently $\theta < \kappa(T)$.
\medskip

\noindent
\underline{Stage B}:

Let $\tau = \tau(T_1) \cup \{P,Q,R,S_{\varphi(\bar x_{[n]},\bar y_{[n]})},
G_n,c_b,Q_\theta,<_\theta,F_i,P_i,F_{1,i}:b \in \bbB,i < \kappa,
\varphi(\bar x_n,\bar y_n) \in \EQ_T\}$, see Definition \ref{b13}(5)
on $\EQ_T$; where the union is without repetitions, 
$P_i,Q_\theta$ unary predicates, $c_b$ an individual 
constant, $R$ binary predicate, $S_{\varphi(\bar x_{[n]})}$ an 
$n$-place function for $\varphi(\bar x_{[n]}) \in \bbL(\tau_T),F_i$ 
unary function for $i < \kappa;F_{1,n}$ is an $n$-place function
symbol, $G_n$ an $n$-place function symbol.

For awhile fix $M_1 \in \Mod_{\bold t}$, note that by \ref{b18}
\mn
\begin{enumerate}
\item[$(*)_1$]  $\|M_1\| =\|M_1\|^{< \kappa} \ge \lambda(T)$.
\end{enumerate}
\mn
Let $M=M_1 \rest \tau(T)$ and let $\cM[M_1]$ be the set $N$ of such
that (for use in other places in $(*)_2$ we do not use ``$\kappa > \aleph_0$"):
\mn
\begin{enumerate}
\item[$(*)_2$]  
\begin{enumerate}
\item[(a)]  $N$ is a $\tau$-expansion of $M_1$
\sn
\item[(b)]  $P^N,Q^N,R,\langle c^N_b:b \in \bbB\rangle$
code $\bbB_T$ and $\uf(\bbB_T)$, see \ref{a4}(3) and 
 e.g. \ref{b24}$(B)(b)(\gamma)-(\eta)$ 
\sn
\item[(c)]
\begin{enumerate}
\item[$(\alpha)$]  $S^N_{\varphi(\bar x_{[m]})}
(\bar a) = \{c^N_b\}$ when $M \models \varphi[\bar a]$; essentially
this says $b = \varphi_b(x_{[m]})/\equiv_T$ for $b \in \bbB_{T,m}$
\sn
\item[$(\beta)$]  $Q^N = \{d_D:D \in \uf(\bbB_T)\}$
  and $R^N = \{(c^N_b,d_D):b \in \bbB$ and $D \in \uf(\bbB),b \in D\}$
\end{enumerate}
\sn
\item[(d)]  for every $i < \kappa$ and $\bar a =
\langle a_j:j <i\rangle \in {}^i M$ for some $b \in N$ we have
$(\forall j < i)(F^N_j(b)=a_j)$ and $b \in P^N_i$
\sn
\item[(e)]  $\langle P^N_i:i < \lambda\rangle$ is a partition of $N$
\sn
\item[(f)]
\begin{enumerate}
\item[$(\alpha)$]  $F^N_{1,m}$ is a function from
${}^m M$ to $Q^N$ such that if $\bar a \in {}^m M$
\then \, $d = F^N_{1,m}(\bar a)$ is the member of $Q^N$
coding $\tp(\bar a,\emptyset,M)$, i.e. 
\sn
\begin{itemize}
\item if $D \in \uf(\bbB_T)$, then we have that $F_{1,m}(\bar a) = d_D$
  if and only if $\tp(\bar a,\emptyset,M) = D$
\end{itemize}
\sn
\item[$(\beta)$]  if $D \in \uf(\bbB_{T,m})$ then
  for some $\bar a \in {}^m M,F^N_{1,m}(\bar a) = d_D$, (recall
 $\bbB_{T,m} \subseteq \bbB_T$)
\end{enumerate}
\sn
\item[(g)]   for any regular $\theta < \kappa_r(T)$ we have:
\sn
\begin{enumerate}
\item[$(\alpha)$]  $Q^N_\theta = \cup\{P^N_i:i \le
  \theta\}$ and $(Q^N_\theta,<^N_\theta)$ is a
  partial order which is a tree with $\theta$ levels isomorphic to 
$({}^{\theta \ge}\|M_1\|,\triangleleft)$ say
$\pi_\theta:{}^{\theta >}\|M_1\| \rightarrow Q^N_\theta$
is such an isomorphism
\sn
\item[$(\beta)$]   let $\bar a^\theta_\eta = \langle
F^N_i(\pi_\theta(\eta)):\ell < \ell g(\bar
  y_{\theta,i})\rangle$ for $\eta \in {}^{\theta \ge}\|M_1\|$
\sn
\item[$(\gamma)$]  $b_1 <^N_\theta b_2$ \Iff \, for some
  $i_1 <i_2 < \theta$ we have $b_1 \in P^N_{i_1},b_2 \in P^N_{i_2}$
and $j < \ell_1 \Rightarrow F^N_j(b_1) = F^N_j(b_2)$
\sn
\item[$(\delta)$]   if $i < \theta,\eta \in {}^i \|M_1\|$
  and $\alpha < \beta < \|M_1\|$ then $N \models \neg(\exists x)
\big((\varphi_{\theta,i}(x,\bar a^\theta_{\eta \char 94 
\langle \alpha \rangle}) \wedge 
\varphi_i(x,\bar a^\theta_{\eta \char 94 \langle \beta\rangle})\big)$
\sn
\item[$(\varepsilon)$]   if $n < \omega,i_0 < \ldots <
  i_{n-1} < \theta,\eta_k \in {}^{(i_k)}\|M_\ell\|$ for $k < n$ and $\eta_0
\triangleleft \eta_1 \triangleleft \ldots \eta_{n-1}$ then
$N \models (\exists x)
(\bigwedge\limits_{k<n} \varphi_{i_k}(x,\bar a^\theta_{\eta_k}))$
\sn
\item[$(\zeta)$]  $F_{\theta,j,i}(\pi(\eta)) = \pi(\eta
  \rest i)$ when $i < j \le \theta,\eta \in {}^j\|M_1\|$
\sn
\item[$(\theta)$]   for every $c \in Q^N_\theta,F^N_\theta(c)$ is
$\pi_\theta(\eta)$ for some $\eta \in {}^{\theta \ge}\|M_1\|$
letting $j_\eta = \ell g(\eta)$ we have
\sn
\begin{itemize}
\item   if $i < j_\eta$ then $N \models
  \varphi_{\theta,i}[c,\bar a^0_{\eta \rest (i_1)}]$
\sn
\item  if $j_\eta < \theta$ then
  $\alpha < \|M_1\| \Rightarrow N \models \neg \varphi_{J_\eta}[c,\bar
  a_{\eta \char 94 \langle \alpha \rangle}]$
\end{itemize}
\sn
\item[$(\iota)$]  $F^N_{\theta,2}$ is a binary function
  such that: if $\eta \in {}^{\theta >}\|M_1\|$ then $\langle
F^N_{\theta,i}(c,\pi_\theta(\eta)):c \in \|M_1\|\rangle$ list with no
  repetitions $\langle \pi_\theta(\eta \char 94 \langle
  \alpha\rangle):\alpha < \|M_1\|\rangle$
\sn
\item[$(\kappa)$]  $F^N_{i,1,\theta}$ or $F^N_{\theta,1}$ is a unary function 
 for every $c \in M,F_{1,\theta}(c)$ is
\sn
\begin{itemize}
\item  $\pi(\eta)$ for some $\eta \in
  {}^{\theta \ge}\|M_1\|$ and for any $i \le \theta,\nu \in
  {}^i\|M_1\|$ we have $c$ realize $\{\varphi_j(x,\bar a^\theta_{\nu
    \rest j}):j<i$ iff $\nu \trianglelefteq \eta\}$
\end{itemize}
\end{enumerate}
\sn
\item[(h)]
\begin{enumerate}
\item[$(\alpha)$]   if $j < \kappa$ has cofinality
  $\theta$, then we have witnesses for clause (d), i.e. if it holds for 
every $j_1 < j$ then it holds for $j$; that is, choose 
$\langle i_j(\iota):\iota < \theta\rangle$, an increasing with limit $j$ and 
 demand:
\newline
\Iff \, $b_i \in M_2$ for $i<j,d \in N$ and 
$F^N_{\theta,2}(d) \in P^N_\theta$ and $\iota < \theta \wedge i_*
< i_j(\iota) \Rightarrow F^N_{i_*}(F^N_\iota(d)) = b_{i_*}$ \then \,
there is $d' \in P_j$ such that $i_* < j \Rightarrow F_{i_*}(d') = b_{i_*}$
\end{enumerate}
\sn
\item[(i)]
\begin{enumerate}
\item[$(\alpha)$]   if $\kappa > \aleph_0$ and 
$\langle a_n:n < \omega\rangle$ is an indiscernible set in $M$ 
\then \, for\footnote{note that when $\kappa > \aleph_0$ we
  can use $G$ a two-place function symbol} some $b,a \mapsto
G^N_2(a,b)$ is a one-to-one function from $M$ onto an indiscernible set which
includes $\{a_n:n < \omega\}$
\sn
\item[$(\beta)$]   if $\kappa = \aleph_0,\bar
  c \in {}^n M,b \in M$ is not algebraic over $\bar c$, then
\sn
\begin{itemize}
\item   $a \mapsto G^N_{n+2}(a,b,\bar c)$ is one-to-one
\sn
\item  $G^N_{n+2}(b,b,\bar c) = \bar b$
\sn
\item   $\{G^N_{n+2}(a,b,\bar c):a \in M\}$ is
an indiscernible set over $\bar c$ based on $\bar c$, all in $M$.
\end{itemize}
\end{enumerate}
\end{enumerate}
\end{enumerate}
\mn
Let $\psi \in \bbL_{\lambda^+,\kappa}[\bbB](\tau)$ be such that:
\mn
\begin{enumerate}
\item[$(*)_3$]  a $\tau$-model $N$ satisfies $\psi$ \Iff \,: for a
  relevant large enough subset $\Lambda$ of
  $\bbL_{\lambda^+,\kappa}[\bbB](\tau)$ of cardinality $\le
  \lambda,\psi = \wedge\{\varphi \in \Lambda$: if $M_1 \in
\Mod_{\bold t}$ and $N \in \cM[M_1]$ then $N \models \varphi\}$; we
may alternatively demand $\psi$ is such that clauses (a)-(h) below hold:
\sn
\begin{enumerate}
\item[(a)]  $N \rest \tau_T$ is a model of $T$, moreover
\sn
\item[(b)]  $N \rest \tau_{T_1}$ is a model of $T_1$
\sn
\item[(c)]  $N \rest \tau_{T_1}$ omits $p$
\sn
\item[(d)]  (e),(f) \quad the parallel of those clauses in $(*)_2$
\sn
\item[(g)] for every $m$, every $m$-type coded by some
  $a \in \bbB_{T,m}$ if $b \in P^N_{2i}$ code 
$\langle a_j:j < 2i\rangle$ satisfies 
$\langle a_{2j},a_{2j+1}:j<i\rangle$ is a $\tau$-elementary 
mapping and $a_{2i} \in N$ then for some $b' \in P_{2i+1}$ and $a_{2i+1}$ 
the element $b'$ code the $\tau$-elementary mapping 
$\langle (a_{2j},a_{2j+1}):j \le i\rangle$
\sn
\item[(h)]  recalling $\kappa > \aleph_0$ if $\langle a_n:n
< \omega\rangle$ is an indiscernible set \then \, for some 
$b,a \mapsto G^N_2(a,b)$ is a one-to-one function from $N$
onto an indiscernible set which includes $\{a_n:n < \omega\}$.
\end{enumerate}
\end{enumerate}
\mn
Now
\mn
\begin{enumerate}
\item[$(*)_4$]
\begin{enumerate}
\item[(a)]  $\psi \in \bbL_{\lambda^+,\kappa}[\bbB]$ indeed
\sn
\item[(b)]   every $M_1 \in \Mod_{\bold t}$ can be
expanded to a model for $\Mod^*_\psi$ (see Definition \ref{b9}(2); this is more
than being a model of $\psi!$)
\sn
\item[(c)]  if $N \in \Mod_\psi$ then $N \rest \tau(T_1) \in \Mod_{\bold t}$.
\end{enumerate}
\end{enumerate}
\mn
[Why?  For clause (a) read $(*)_3$.  For clause (b) read $(*)_2 +
(*)_3$.  For clause (c), first why $M_1 = N \rest \tau_{T_1}$ is a
model of $T_1$?   Since $M_1 \in \Mod_t$ and $N \in \cM[M_1]$, we have
that $N \upharpoonleft \tau(T_1)$ is a $\tau$-expansion of $M_1$ by $(*)_2(a)$.
Second, why $M_1$ omit
$p_{\bold t}$?  Recalling $(*)_2(f)(\alpha) + (\beta)$ and choice of
$\psi$ this should be clear.  Third, why is $M = N \rest \tau_T$ saturated?
  It realizes every $p \in D_m(T) = \bold S^m(\emptyset,M)$, by
$(*)_2(f)$, it is $\kappa$-sequence-homogeneous by $(*)_3(g)$ hence
 is $\kappa$-saturated.  By $(*)_3(h)$, every indiscernible subset
 $\bold I$ of cardinal $\aleph_0$ 
can be extended to one of cardinality $\|M\|$.  By 
the last two sentences, $M$ is saturated by Case 1 of \ref{a9}.]

So we are done.
\end{PROOF}

\begin{claim}
\label{d8}
Like \ref{d4}, but $T$ is superstable and $\lambda(T) \le \lambda$.
\end{claim}

\begin{PROOF}{\ref{d8}}
Here the proof ``why $M=N \rest \tau_T$ is saturated inside the proof
of $(*)_4(c)$ is different.
There is a saturated $M_* \in \Mod_T$ of cardinality $\le \lambda$ and
we can demand on $\psi$ that $N \models \psi$ implies $M_*$ is elementarily
embeddable into $N \rest \tau_T$ and $N \rest \tau_T$ is
$\aleph_0$-sequence homogeneous.  

Note that
\mn
\begin{enumerate}
\item[$(*)$]  if $M_* \prec M \in \Mod_T$ and $M$ is
  $\aleph_0$-sequence homogeneous implies $M$ is
  $\aleph_\varepsilon$-saturated, see \ref{b13}(0).
\end{enumerate}
\mn
In this case $(*)_2(i)(\beta)$ of the proof of \ref{d4} implies $M$ is
saturated because by case 2 of \ref{a9}
\mn
\begin{enumerate}
\item[$(*)$]  $M$ is saturated \when \,: if $M$ is
  $\aleph_\varepsilon$-saturated and for every finite $A \subseteq M$
  and $a \in M \backslash \acl(A)$ there is an indiscernible set $\cI
  \subseteq M$ over $A$ of cardinal $\|M\|$ based on $A$
  (i.e. $\Av(M,\bold I)$ does not fork over $A$) to which $a$ belongs.
\end{enumerate} 
\end{PROOF}

\begin{claim}
\label{d12}
1) Like \ref{d4} but $T$ is superstable and $2^{\aleph_0} \le
   \lambda$.

\noindent
2) Like \ref{d4}, but $T$ superstable and $|D(T)| > |T|$.
\end{claim}

\begin{PROOF}{\ref{d12}}
As the proof of \ref{d8} the problem is how $\psi$ guarantees ``$N
\rest \tau_T$ is $\aleph_\varepsilon$-saturated".  As the model is
$\aleph_0$-sequence homogeneous it suffices
\mn
\begin{enumerate}
\item[$(*)$]  for every $m$ and $D \in \uf(\bbB_{T,m+1})$ equivalently
  $p \in D_{m+1}(T)$ for some
  $\bar a \char 94 \langle c \rangle \in {}^{m+1}N$ realizing $p$, we
  have: if $N \rest \tau_T \prec M'$ and 
$c' \in M'$ realizes $\tp(c,\bar a,N \rest \tau_T)$
  \then \, some $c'' \in N \rest \tau_T$ realizes $\stp(c',\bar a,M')$
  in $M'$.
\end{enumerate}
\mn
Let $p = \tp(c,\bar a,M)$ and we let $\lambda_* = \lambda(p),\langle
E_\alpha(x_0,x_1;\bar y_{[m]}):\alpha < \lambda_*\rangle$, see
\cite[Ch.III,5.1,pg.123]{Sh:c}.
\medskip

\noindent
\underline{Case 1}:  $\lambda_* = \aleph_0$

If $2^{\aleph_0} \le \lambda$ this is easy.  If $|D(T)| > |T|$ then
for some $m$ there is an independent sequence $\langle \varphi_n(\bar
x_{[m]}):n < \omega\rangle$ of formulas of $\bbL(\tau_T)$ over $T$;
(that is, if $M \in \Mod_T$ then any non-trivial finite Boolean
combination of them is realized in $M$) and
we continue as in the second case.
\medskip

\noindent
\underline{Case 2}:  $\lambda_* > \aleph_0$

In this case by \cite[Ch.III,5.9,5.10,pg.126]{Sh:c} there is
a sequence of length $\lambda_*$ of formulas of the form $\varphi[x,\bar a]$ 
independent in $\gC_T$.  Hence
there is an independent over $T$ sequence $\langle \varphi_i(x,\bar
y_{[m]}):i < \lambda_*\rangle$ of formulas from $\bbL(\tau_T)$, so
$\bbB^{\fr}_{\lambda_*}$ is embeddable into $\bbB_{T,m+1}$.  So $\psi$
says that the Boolean Algebra $\cP(\lambda_*)$ is interpreted in $N$
for every relevant $\lambda_*$, but $\lambda_* \le |T|$.

\relax From this it is easy to have $\psi$ ensuring $(*)$.
\end{PROOF}
\bigskip

\subsection {Coding $\psi \in \bbL_{\lambda^+,\kappa}[\bbB_T]$} \
\bigskip

\begin{hypothesis}  
\label{e2}
\mn
\begin{enumerate}
\item[$(a)$]  $T$ is a complete first order theory,
\sn
\item[$(b)$]  $\lambda \ge |T|,\lambda^+ \ge \kappa$
\sn
\item[$(c)$]  $\bbB = \bbB_T$.
\end{enumerate}
\end{hypothesis}

\begin{claim}  
\label{e4}
Assume $\psi \in \bbL_{\lambda^+,\kappa}[\bbB]$ and 
$\kappa = \kappa_r(T) < \infty$ so $T$ is stable.

There is $\bold t = (T,T_1,p) \in \bold N_{\lambda,T}$ such that
$\tau(T_1) \supseteq \tau(\psi)$ and $\Mod_{\bold t} = \{N \rest
\tau(\psi):N \in \Mod_\psi[\bbB]\}$.
\end{claim}

\begin{PROOF}{\ref{e4}}  
We apply \ref{b24} to $\bbB$ and $\psi$ and get $(\tau_1,T_1,p(*),\bar
F,\bar P)$ as in \ref{b20}, \ref{b24} and \wilog \, $\tau_1 \cap \tau(T) =
\emptyset$.  Now we immitate the proof of \ref{d4}.  [Referee 2.4]
\end{PROOF}
\bigskip

\subsection {Elaborating Case C}\
\bigskip

\noindent
In \S(2B) we treat most theories $T$ but not all.
The remaining case is
\begin{hypothesis}
\label{f2}
\mn
\begin{enumerate}
\item[$\boxplus$]  $(a) \quad T$ is superstable of cardinality $\lambda$
\sn
\item[${{}}$]  $(b) \quad \lambda(T) > \lambda$
\sn
\item[${{}}$]  $(c) \quad 2^{\aleph_0} > \lambda$
\sn
\item[${{}}$]  $(d) \quad \lambda \ge |D(T)|$.
\end{enumerate}
\end{hypothesis}

\begin{claim}
\label{f6}
There are $m,M \in \Mod_T$ and $\bar a \in {}^m M$ such that
$\{\stp(c,\bar a,M):c \in M\}$ is of cardinality $2^{\aleph_0}$.
\end{claim}

\begin{PROOF}{\ref{f6}}
Should be clear. [Referee 2.16]
\end{PROOF}

\begin{definition}
\label{f10}
For any model $M$ and a sequence $\bar a$ from $M$ (or a set
$\subseteq$), let $\bbB_{M,\bar a,m}$ be the Boolean Algebra of
subsets of ${}^m M$ of the form $\varphi(M,\bar c)$, where
$\varphi(\bar x_{[m]},\bar z) \in \bbL(\tau_M),\bar b \in {}^{\ell
  g(\bar z)}M$ and $\varphi(\bar x,\bar c)$ is almost over $\bar a$
which means: for some $\vartheta(\bar x_{[m]},\bar y_{[m]},\bar z) \in
\bbL(\tau_M)$ we have:
\mn
\begin{itemize} 
\item   in $M,\vartheta(\bar x_{[m]},\bar y_{[m]},\bar a)
\vdash \varphi(\bar x_{[m]},\bar c) \equiv \varphi(\bar y_{[m]},\bar c)$
\sn
\item   $\vartheta(\bar x_{[m]},\bar y_{[m]},\bar a)$
  defines in $M$ an equivalence relation with finitely many
  equivalence classes.
\end{itemize}
\end{definition}

\begin{claim}
\label{f12}
For $T$ as in \ref{f2}, letting $M,\bar a,m$ be as in \ref{f6} and
$\bbB = \bbB_{M,\bar a,m}$ the result of \ref{d4} and Theorem \ref{c2}
hold if we use $\bbB$ instead of $\bbB_T$.
\end{claim}

\begin{PROOF}{\ref{f12}}
As above, really $m=1$ suffice; in particular if $p \in \bold S(\bar
a,M),\bar a \in {}^m M,M \in \Mod_T$ then $\lambda_*(p) \le \aleph_0$
(otherwise by Lemma 5.9, 5.10 and 5.11 \cite[Ch.III]{Sh:c} we have
$|\bold S^{2m}(\bar a,m)| \ge 2^{\lambda_*(p)} > \lambda$, contradiction).
\end{PROOF}
\newpage

\bibliographystyle{alphacolon}
\bibliography{lista,listb,listc,listd,liste,listf,listv,listx,listy,listz}

\end{document}